\documentclass{article}

\usepackage[a4paper,top=2.54cm,bottom=2.54cm,left=3.17cm,right=3.17cm,%
            includehead,includefoot]{geometry}

\usepackage{amsmath,amssymb,amsfonts,amsthm}
\usepackage{graphicx}
\usepackage{subfigure}
\usepackage{float}
\usepackage{xcolor}
\usepackage[numbers,square,sort&compress]{natbib}
\usepackage{hyperref}
\hypersetup{colorlinks,citecolor=blue,linkcolor=blue,breaklinks=true}
\usepackage{booktabs}
\usepackage{colortbl}
\usepackage{caption}
\usepackage{enumitem}
\usepackage{algorithm}
\usepackage{algpseudocode}
\usepackage{array}
\usepackage{bbding}
\usepackage{fancyhdr} 
\usepackage{fancyvrb} 
\usepackage{longtable}
\usepackage{listings}
\usepackage{rotating,rotfloat} 
\usepackage{yhmath} 
\usepackage{tikz}
\usetikzlibrary{3d, positioning, shapes.geometric, arrows.meta, calc}
\usepackage{stmaryrd}
\usepackage{multirow, bigstrut}
\usepackage{graphicx}  
\usepackage{tikz}
\usetikzlibrary{positioning}
\usepackage{subcaption}
\usepackage{adjustbox} 
\usepackage{arydshln}   

\usepackage{authblk}  

\usepackage{cleveref}

\usetikzlibrary{positioning}
\tikzset{
  input/.style={draw, rectangle, fill=blue!10, rounded corners, minimum height=1.2em, minimum width=4em},
  nn/.style={draw, rectangle, fill=green!10, rounded corners, minimum height=1.5em, minimum width=4em},
  operator/.style={draw, circle, fill=gray!10, minimum size=2em},
  output/.style={draw, rectangle, fill=orange!20, rounded corners, minimum height=1.2em, minimum width=4em}
}

\usepackage{fancyhdr}
\usepackage{booktabs}  
\allowdisplaybreaks

\newtheorem{definition}{Definition}[section]

\newtheorem{remark}{Remark}

\newcommand{\bbm}{\begin{bmatrix}}
\newcommand{\ebm}{\end{bmatrix}}

\title{Neural Evolutionary Kernel Method: A Knowledge-Guided Framework for Solving Evolutionary PDEs}

\author[1]{Shuo Ling\thanks{\href{mailto:lingshuo1@sjtu.edu.cn}{lingshuo1@sjtu.edu.cn}}}
\author[2]{Wenjun Ying\thanks{\href{mailto:wying@sjtu.edu.cn}{wying@sjtu.edu.cn}}}
\author[3]{Zhen Zhang\thanks{Corresponding author. \href{mailto:zhangz@sustech.edu.cn}{zhangz@sustech.edu.cn}}}
\affil[1]{School of Mathematical Sciences, Shanghai Jiao Tong University, Minhang, Shanghai 200240, PR China}
\affil[2]{School of Mathematical Sciences, MOE-LSC and Institute of Natural Sciences, Shanghai Jiao Tong University, Minhang, Shanghai 200240, PR China}
\affil[3]{Department of Mathematics, Southern University of Science and Technology, Nanshan, Shenzhen 518000, Guangdong, PR China}

\date{\today}

\begin{document}
\maketitle

\begin{abstract}
Numerical solution of partial differential equations (PDEs) plays a vital role in various fields of science and engineering. In recent years, deep neural networks (DNNs) have emerged as a powerful tool for solving PDEs, leveraging their approximation capabilities to handle complex domains and high-dimensional problems. Among these, operator learning has gained increasing attention by learning mappings between function spaces using DNNs.
This paper proposes a novel approach, termed the Neural Evolutionary Kernel Method (NEKM), for solving a class of time-dependent partial differential equations (PDEs) via deep neural network (DNN)-based kernel representations. By integrating boundary integral techniques with operator learning, prior mathematical information of time-dependent partial differential equations (PDEs) is embedded into the design of neural network architectures for predicting their solutions, enhancing both computational efficiency and solution accuracy. Numerical experiments on the heat, wave, and Schr\"{o}dinger equations demonstrate that the Neural Evolutionary Kernel Method (NEKM) achieves high accuracy and favorable computational efficiency. Furthermore, the operator learning framework inherently supports the simultaneous prediction of solutions to multiple PDEs with different coefficients, rendering its capability for solving random PDEs. 

\end{abstract}

\section{Introduction} \label{Introduction}
Partial differential equations (PDEs) play a central role in modeling a wide range of physical and chemical phenomena, including diffusion, heat transfer, and fluid dynamics. In most practical settings, exact solutions of PDEs are unavailable, and numerical approximations serve as indispensable tools for scientific computation. Over the past century, numerous numerical methods have been developed for solving PDEs, such as finite difference, finite element, finite volume, and spectral methods. Despite their success, these classical approaches often suffer from the curse of dimensionality \citep{weinan2019barron}, and the construction of high-order, stable numerical schemes becomes increasingly challenging in complex geometries.

Motivated by the Universal Approximation Theorem \cite{cybenko1989approximation}, neural networks have been widely explored for function approximation \citep{weinan2019barron, barron1993universal, shen2022deep, zhang2022neural}, leading to the rapid development of deep neural network (DNN)-based methods for solving PDEs. Early attempts can be traced back to \cite{dissanayake1994neural}, while the deep learning era has significantly accelerated progress in this field. Representative approaches include the deep Galerkin method (DGM) \citep{sirignano2018dgm} and physics-informed neural networks (PINNs) \citep{raissi2019physics}, which formulate PDE solving as the minimization of equation residuals sampled at randomly selected spatial points. Boundary conditions are enforced either by explicitly modifying network architectures \citep{berg2018unified} or by incorporating penalty terms into the loss function, where the choice of penalty parameters remains problem-dependent. To address the deterioration caused by high-order derivatives in the loss, the deep mixed residual method (MIM) \citep{lyu2022mim} reformulates high-order PDEs as first-order systems. Alternatively, the deep Ritz method \citep{yu2018deep} adopts weak formulations to recast PDEs into variational optimization problems, although this limits its applicability to PDEs with suitable energy structures. Recent works have also provided convergence analyses for PINNs and the deep Ritz method \citep{lu2022machine, CiCP-31-1020, CiCP-31-1272}. Meanwhile, incorporating mathematical prior knowledge into network design has attracted increasing attention, including approaches that embed the Boundary Integral Method (BIM) \citep{atkinson1997numerical} into neural architectures \citep{CSIAM-AM-4-275,lin2023bi}. For time-dependent PDEs, Su et al.~\cite{su2025spikestablephysicsinformedkernel} proposed the Stable Physics-Informed Kernel Evolution (SPIKE) method to resolve shocks in hyperbolic conservation laws. Despite these advances, most DNN-based PDE solvers are tailored to fixed problem instances and lack the ability to efficiently generalize across varying source terms, boundary conditions, or PDE parameters without retraining.

Rather than approximating pointwise solutions, operator learning seeks to learn the global mapping from input functions—such as coefficients, source terms, or boundary data—to solution functions of PDEs. This paradigm enables the approximation of mappings between infinite-dimensional function spaces in a mesh-free manner \cite{lu2021learning, li2021fourier, Bhattacharya2021, doi:10.1137/20M133957X}. Notable examples include the Fourier Neural Operator (FNO) \citep{li2021fourier}, which performs operator learning in the frequency domain using Fourier transforms, and DeepONet \citep{lu2021learning}, which employs a branch–trunk architecture to encode input functions and spatial coordinates. Other approaches, such as DeepGreen \citep{gin2021deepgreen} and MOD-Net \citep{CiCP-32-299}, incorporate Green’s function theory by learning continuous kernel representations of solution operators. While these operator learning methods demonstrate strong generalization across input functions, they are often treated as black-box models and may not fully exploit the underlying mathematical structure of PDEs. As a consequence, such approaches typically require relatively large amounts of training data to achieve satisfactory performance, and the resulting accuracy may still be limited for certain classes of parameterized elliptic problems.
 
The primary objective of this work is to construct an operator learning framework that can be effectively applied to time-dependent PDEs, such as the heat equation, wave equation, and Schrödinger equation. To this end, a knowledge-guided operator learning framework, termed the \emph{Neural Evolutionary Kernel Method} (NEKM), is developed. When solving time-evolution problems, each time step typically requires the solution of a time-discretized PDE. As a result, the corresponding equations may involve varying source terms, boundary conditions, and even different PDE parameters when adaptive or variable time step sizes are employed. This setting naturally calls for an efficient and accurate operator learning approach that, on a fixed computational domain, can accommodate a class of PDEs with varying right-hand sides, boundary data, and equation parameters. To this end, the time-discretized elliptic problem is decomposed into two subproblems, corresponding to the contributions from the source term and the boundary condition, respectively. For the component associated with the source term under homogeneous boundary conditions, a class of neural network architectures is designed based on the structure of volume integral representations. Trained in a supervised manner, this construction enables the accurate approximation of the solution operator associated with volumetric forcing and demonstrates favorable approximation performance. For the component induced by nonhomogeneous boundary conditions, a direct learning of the mapping from boundary data to the solution is not pursued. Such a direct approach would involve a significant mismatch between the input and output dimensions, as the boundary data are defined on $\partial\Omega$ while the solution resides in the full domain $\Omega$, and the corresponding solution operator is only weakly characterized. Instead, the boundary integral method (BIM) offers a natural and mathematically grounded alternative by representing the solution as a boundary integral involving an unknown density function defined solely on the boundary. This density is determined by solving a second-kind boundary integral equation, whose associated solution operator is well-conditioned and free of singularities. Motivated by these observations, a central strategy of this work is to learn the solution operator of the boundary integral equation itself, namely, the mapping from the boundary condition and PDE parameters to the corresponding boundary density. Since both the input boundary data and the output density function are defined on $\partial\Omega$, this approach effectively reduces the dimensionality of the learning problem by one and leverages the favorable analytical properties of second-kind integral operators, thereby enabling high-accuracy approximation. Once the elliptic solution operator is learned with sufficient precision, it can be repeatedly invoked at successive time steps, allowing the accurate numerical solution of time-dependent PDEs through time discretization.

The BIM \cite{atkinson1997numerical, yu2002natural, pozrikidis1992boundary} leverages the second and third Green’s identities to reformulate the PDE as an integral equation posed solely on the boundary of the domain. Once this boundary integral equation is solved, the solution to the original PDE can be reconstructed accordingly. As a competitive alternative to both classical and modern methods, the kernel-free boundary integral (KFBI) method — and its extensions that combine GPU acceleration and neural networks — have demonstrated notable advantages in recent years \cite{ying2007kernel, ying2014kernel, xie2020high, zhao2023kernel, zhou2023kernel, tan2024gpu, ling2025hybrid}. These traditional boundary integral approaches, along with their recent extensions that incorporate deep learning techniques, have offered important conceptual and methodological inspiration for the development of our operator learning framework.

The main contributions of this work are summarized as follows:
\begin{itemize}
  \item A mathematically informed operator learning framework is proposed for linear elliptic PDEs on fixed computational domains, with explicit accommodation of varying source terms, boundary conditions, and PDE parameters. The framework is particularly motivated by time-dependent PDEs, which can be reformulated as sequences of parameterized elliptic problems through implicit temporal discretization.

  \item Exploiting the linearity of elliptic operators, the solution operator is decomposed into two components corresponding to the source term and the boundary condition, respectively. For the source-driven component under homogeneous boundary conditions, a class of neural network architectures inspired by volume integral representations is developed and trained using supervised learning, yielding accurate approximations of the associated solution operator.

  \item For the boundary-driven component, a novel learning strategy based on BIM is introduced. Instead of directly learning the mapping from boundary data to the solution in the full domain, the proposed approach learns the solution operator of a second-kind boundary integral equation, mapping boundary conditions and PDE parameters to the corresponding boundary density. This strategy leverages the favorable analytical properties of second-kind integral operators and reduces the dimensionality of the learning problem, leading to improved accuracy.

  \item By incorporating PDE parameters as additional inputs, a single trained model can be applied to a family of parameterized elliptic operators without retraining. Once trained, the learned elliptic solver can be repeatedly invoked as an evolution kernel for time-dependent PDEs, enabling accurate and efficient numerical simulation over multiple time steps.

  \item The operator-based nature of the proposed framework further facilitates extensions to related tasks such as uncertainty quantification (UQ) \cite{martina2021bayesian} and inverse problems \cite{gao2024adaptive}, where rapid evaluation of the solution operator across varying inputs is essential. Additionally, thanks to the tensor-based design of the \texttt{PyTorch} framework, models can efficiently solve multiple equations in parallel, offering a clear advantage over traditional numerical solvers. 
\end{itemize}

The remainder of this paper is organized as follows. Section~\ref{Preliminaries} reviews the theoretical foundations of the boundary integral method (BIM), which underpin the proposed algorithmic design. Section~\ref{Methodology} details the Neural Evolutionary Kernel Method (NEKM), a unified framework that synergistically combines operator learning with BIM to improve generalization and computational scalability. Section~\ref{Numerical Experiments} presents extensive numerical experiments validating the accuracy, efficiency, and versatility of NEKM. Finally, Section~\ref{Conclusion} discusses the method's strengths, current limitations, and promising avenues for future research.

\section{Preliminaries} \label{Preliminaries}
This section reviews essential aspects of the BIM that underpin the NEKM. For clarity and conciseness, the exposition focuses on the application of BIM to boundary value problems (BVPs) defined on simply connected, bounded domains in two spatial dimensions, subject to Dirichlet boundary conditions. The treatment of Neumann boundary conditions, as well as the extension to three-dimensional domains, follows standard theory and can be found in classical references on the subject, e.g., \cite{atkinson1997numerical}.

Let $\Omega \subset \mathbb{R}^2$ be a simply connected bounded domain. Let $\sigma \equiv \sigma(\mathbf{x})$ be a symmetric and positive definite (SPD) diffusion tensor and $\kappa \equiv \kappa(\mathbf{x})$ be a non-negative reaction coefficient. Then we introduce an (negative) elliptic operator
\begin{equation} \label{eq0606_1}
\mathcal{L} \equiv \nabla \cdot \sigma(\mathbf{x}) \nabla-\kappa(\mathbf{x}) .
\end{equation}
Let $H^{r-1}(\Omega), H^{r+1 / 2}(\partial \Omega)$ be the standard Sobolev spaces, with $r>k$ and $k$ being a positive integer. Given a source function $f \in H^{r-1}(\Omega)$ and boundary data $g^D \in H^{r+1 / 2}(\partial \Omega)$, this elliptic equation is considered
\begin{equation} \label{eq0606_2}
\mathcal{L} u(\mathbf{x})=f(\mathbf{x}) \quad \text { in } \Omega,
\end{equation}
with the pure Dirichlet boundary condition
\begin{equation} \label{eq0606_3}
u(\mathbf{x})=g^D(\mathbf{x}) \quad \text { on } \partial \Omega.
\end{equation}
Eqs. \eqref{eq0606_2} and \eqref{eq0606_3} form an interior Dirichlet boundary value problem. Due to the linearity of $\mathcal{L}$, the solution can be expressed as the sum of $u_1$ and $u_2$, which satisfy problems \eqref{eq0606_4} and \eqref{eq0606_5} respectively.
\begin{equation}\label{eq0606_4}
\begin{cases}
\mathcal{L} u_1 (\mathbf{x}) = f (\mathbf{x}) & \mbox{in}\ \Omega , \\
u_1 (\mathbf{x}) = 0 & \mbox{on}\ \partial \Omega,
\end{cases}
\end{equation}

\begin{equation}\label{eq0606_5}
\begin{cases}
\mathcal{L} u_2 (\mathbf{x}) = 0 & \mbox{in}\ \Omega , \\
u_2 (\mathbf{x}) = g^D (\mathbf{x}) & \mbox{on}\ \partial \Omega.
\end{cases}
\end{equation}

Suppose $G_0(\mathbf{x}, \mathbf{y})$ denotes the fundamental solution of the operator $\mathcal{L}$, and $G(\mathbf{x}, \mathbf{y})$ denotes the Green's function associated with $\mathcal{L}$ in the domain $\Omega$ subject to the Dirichlet boundary conditions, as described in Appendix \ref{appendix A}. Then the solution $u_1$ to equation \eqref{eq0606_4} can be expressed as
\begin{equation} \label{eq0607_1}
u_1(\mathbf{x}) = \int_{\Omega} G(\mathbf{x}, \mathbf{y}) f(\mathbf{y}) d\mathbf{y}, \quad \mathbf{x} \in \Omega.
\end{equation}
The solution $u_2$ to equation \eqref{eq0606_5} can be expressed as a double layer potential
\begin{equation} \label{eq0607_2}
u_2(\mathbf{x}) = \int_{\partial \Omega} \frac{\partial G_0(\mathbf{x}, \mathbf{y})}{\partial n_\mathbf{y}} \varphi(\mathbf{y}) d s_\mathbf{y},\quad \mathbf{x} \in \Omega,
\end{equation}
where $\varphi$ is the solution of the Fredholm boundary integral equation of the second kind \cite{kress1999linear, hsiao2008boundary}:
\begin{equation} \label{eq0607_3}
\frac{1}{2} \varphi(\mathbf{x}) + \int_{\partial \Omega} \frac{\partial G_0(\mathbf{x},\mathbf{y})}{\partial n_\mathbf{y}} \varphi(\mathbf{y}) d s_\mathbf{y} = g^D(\mathbf{x}),\quad \forall \mathbf{x} \in \partial \Omega.
\end{equation}
\begin{remark}
The fundamental solution \( G_0(\mathbf{x}, \mathbf{y}) \) of the operator \( \mathcal{L} \) can usually be derived analytically. However, for general domains \( \Omega \), it is often extremely difficult to obtain an explicit expression for the corresponding Green's function \( G(\mathbf{x}, \mathbf{y}) \). Consequently, in practical computations, it is rare to directly use formula \eqref{eq0607_1} to evaluate \( u_1 \). Nevertheless, formula \eqref{eq0607_1} provides valuable structural insights, which can be leveraged to design neural networks as discussed below.
\end{remark}

\section{Methodology} \label{Methodology}
This section presents the implementation details of the NEKM. An operator learning framework is first developed for the interior Dirichlet boundary value problem given by \eqref{eq0606_2}--\eqref{eq0606_3}, enabling the trained model to map a source term $f$ and Dirichlet boundary data $g^D$ to the corresponding PDE solution with high accuracy and computational efficiency. Subsequently, the extension of the trained model to time-dependent problems is described, thereby completing the full NEKM pipeline.

\subsection{Operator Learning} \label{Operator Learning}
To predict the solution of the interior Dirichlet boundary value problem \eqref{eq0606_2}--\eqref{eq0606_3}, operator learning is performed separately for equations \eqref{eq0606_4} and \eqref{eq0606_5}. Specifically, for equation \eqref{eq0606_4}, a neural network is designed to map the source term $f$ to $u_1$. For equation \eqref{eq0606_5}, the boundary integral equation \eqref{eq0607_3} serves as the starting point for network design. It is worth noting that, for a general symmetric positive-definite (SPD) diffusion tensor $\sigma$ and non-negative reaction coefficient $\kappa$, the integral equation \eqref{eq0607_3} is non-singular \cite{ying2007kernel}. This property allows the construction of a mapping from the Dirichlet boundary condition $g^D$ to the solution $\varphi$, thereby approximating the associated solution operator. Once this mapping is learned by a neural network, the solution to \eqref{eq0606_5} can be recovered using the representation \eqref{eq0607_2}. In the following, a detailed exposition of the operator learning strategies for these two components is provided.

\subsubsection{Operator Learning for Source Terms} \label{Operator Learning for Source Terms}
To approximate this integral \eqref{eq0607_1} numerically, consider a `quadrature rule' with $n$ quadrature nodes:
\begin{equation} \label{eq0607_4}
\int_{\Omega} G(\mathbf{x}, \mathbf{y}) f(\mathbf{y}) \, d\mathbf{y} \approx \sum_{i=1}^{n} \omega_i G(\mathbf{x}, \mathbf{y}_i) f(\mathbf{y}_i),
\end{equation}
where $\{\omega_i\}_{i=1}^n$ and $\{\mathbf{y}_i\}_{i=1}^n$ are the quadrature weights and nodes in the domain $\Omega$. For brevity, denote
\[
\widetilde{G}_i^{\mathbf{x}} \triangleq \omega_i G(\mathbf{x}, \mathbf{y}_i), \quad \widetilde{f}_i \triangleq f(\mathbf{y}_i).
\]
Let
\[
\widetilde{G}(\mathbf{x}) \triangleq 
\begin{pmatrix}
\widetilde{G}_1^{\mathbf{x}} \\ \widetilde{G}_2^{\mathbf{x}} \\ \vdots \\ \widetilde{G}_n^{\mathbf{x}}
\end{pmatrix}
\in \mathbb{R}^n, \qquad
\widetilde{f} \triangleq 
\begin{pmatrix}
\widetilde{f}_1 \\ \widetilde{f}_2 \\ \vdots \\ \widetilde{f}_n
\end{pmatrix}
\in \mathbb{R}^n,
\]
then the approximated solution becomes a dot product:
\begin{equation} \label{eq0607_5}
u_1(\mathbf{x}) \approx \widetilde{G}(\mathbf{x})^\top \cdot \widetilde{f}, \quad \mathbf{x} \in \Omega.
\end{equation}

However, for general domains $\Omega$, the analytical form of the Green's function $G$ is often unavailable, making it difficult to obtain $\widetilde{G}(\mathbf{x})$ explicitly. To address this, the mapping from $f$ to $u_1$ is learned via a neural network designed based on the above formulation. Suppose $M$ input source functions $f_1, f_2, \ldots, f_M$ are given. Since functions cannot be directly fed into neural networks, $N$ evaluation points $\{\mathbf{x}_j\}_{j=1}^N \subset \Omega$ are sampled, and each function $f_i$ is encoded by its values at these points:
\[
\texttt{f}_{i,j} = f_i(\mathbf{x}_j), \quad \texttt{f} \in \mathbb{R}^{M \times N}.
\]
Similarly, define the output matrix $\texttt{u} \in \mathbb{R}^{M \times N}$ by
\[
\texttt{u}_{i,j} = \hat{u}_i(\mathbf{x}_j),
\]
where $\hat{u}_i$ denotes the predicted solution corresponding to source $f_i$ evaluated at $\mathbf{x}_j$. In addition, the sampling coordinates are collected into a matrix
\[
\texttt{x} = 
\begin{pmatrix}
\mathbf{x}_1^\top \\
\mathbf{x}_2^\top \\
\vdots \\
\mathbf{x}_N^\top
\end{pmatrix}
\in \mathbb{R}^{N \times d},
\]
where $d = 2$ is the spatial dimension of the domain $\Omega$. To approximate the quadrature-based representation described above, two sub-networks are constructed. A network $\mathcal{NN}_f$ maps the input $\texttt{f} \in \mathbb{R}^{M \times N}$ to a latent representation
  \[
  \mathcal{NN}_f(\texttt{f}) = \begin{bmatrix}
  \widetilde{f}_1^\top \\
  \widetilde{f}_2^\top \\
  \vdots \\
  \widetilde{f}_M^\top
  \end{bmatrix}
  \in \mathbb{R}^{M \times n}.
  \]
And another network $\mathcal{NN}_G$ maps the spatial coordinates $\texttt{x} \in \mathbb{R}^{N \times d}$ to the Green's functions at many evaluation points
  \[
  \mathcal{NN}_G(\texttt{x}) = \begin{bmatrix}
  \widetilde{G}(\mathbf{x}_1)^\top \\
  \widetilde{G}(\mathbf{x}_2)^\top \\
  \vdots \\
  \widetilde{G}(\mathbf{x}_N)^\top
  \end{bmatrix}
  \in \mathbb{R}^{N \times n}.
  \]
Finally, the complete neural network is defined as
\begin{equation} \label{eq0607_6}
\mathcal{NN}_1(\texttt{f}, \texttt{x}) := \mathcal{NN}_f(\texttt{f}) \cdot \mathcal{NN}_G(\texttt{x})^\top \in \mathbb{R}^{M \times N},
\end{equation}
where the $(i,j)$-th entry corresponds to the inner product $\widetilde{f}_i^\top \cdot \widetilde{G}(\mathbf{x}_j)$, representing the predicted solution $\hat{u}_i(\mathbf{x}_j)$ for $f_i$ evaluated at location $\mathbf{x}_j$.

\begin{remark}
Strictly speaking, for any $\mathbf{x} \in \Omega$, the integral in \eqref{eq0607_1} exhibits a singularity at $\mathbf{y} = \mathbf{x}$ (if $f$ is continuous on $\bar{\Omega}$) due to the nature of the Green's function. This implies that the singularity location depends on the evaluation point $\mathbf{x}$, and therefore it is difficult to define a universally valid and convergent quadrature rule such as \eqref{eq0607_4} for all $\mathbf{x} \in \Omega$. This is also the reason we refer to it as a `quadrature rule' rather than a strict numerical integration method. Nevertheless, for the purpose of constructing our neural network architecture, the formulation in \eqref{eq0607_4} is sufficiently insightful and has led to promising results in our numerical experiments.
\end{remark}

\begin{figure}[htbp]
\centering
\begin{tikzpicture}[node distance=1.8cm and 2cm,>=latex]
  \node[input] (f_in) at (0,0) {$\texttt{f} \in \mathbb{R}^{M \times N}$};
  \node[input, below=of f_in] (k_in) {$\texttt{k} \in \mathbb{R}^{M \times n_p}$};
  \node[input, below=of k_in] (x_in) {$\texttt{x} \in \mathbb{R}^{N \times d}$};

  \node[nn, right=of f_in] (NNf) {$\mathcal{NN}_f$};
  \node[nn, right=of k_in] (NNk) {$\mathcal{NN}_k$};
  \node[nn, right=of x_in] (NNx) {$\mathcal{NN}_G$};

  \node[operator, right=1.8cm of NNk] (otimes) {$\odot$};
  \node[operator, right=2cm of otimes] (dot) {$\cdot$};

  \node[output, right=of dot] (out) {$\hat{u} \in \mathbb{R}^{M \times N}$};

  \draw[->] (f_in) -- (NNf);
  \draw[->] (k_in) -- (NNk);
  \draw[->] (x_in) -- (NNx);

  \draw[->] (NNf.east) to[out=0,in=135] (otimes.north west);
  \draw[->] (NNk.east) -- (otimes.west);
  \draw[->] (otimes.east) -- (dot.west);
  \draw[->] (NNx.east) to[out=0,in=225] (dot.south west);
  \draw[->] (dot.east) -- (out.west);
\end{tikzpicture}
\caption{Schematic diagram of the network architecture described in \eqref{eq0608_1}. In this framework, the sub-networks $\mathcal{NN}_f$, $\mathcal{NN}_k$, and $\mathcal{NN}_G$ can be instantiated using architectures such as multi-layer perceptrons (MLPs), residual neural networks (ResNets), or other function-approximating models.}
\label{fig0608_1}
\end{figure}
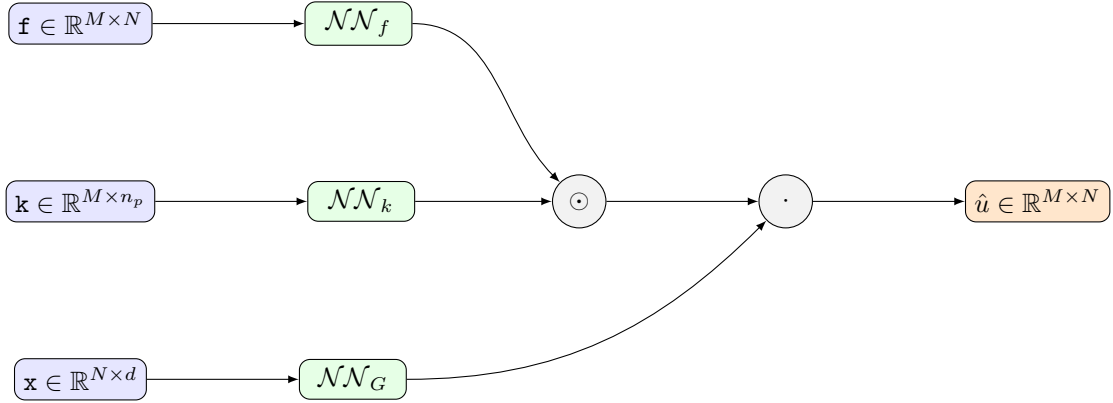

Furthermore, we aim to incorporate certain parameters in the operator $\mathcal{L}$ as part of the network input. For instance, when $\mathcal{L} = \Delta - k \mathcal{I}$ with $k > 0$, where $\mathcal{I}$ denotes the identity operator, $k$ is included as an additional input to the network. Building upon the previous network structure \eqref{eq0607_6}, the architecture is extended to accommodate such parameters:
\begin{equation} \label{eq0608_1}
\mathcal{NN}_1(\texttt{f}, \texttt{k}, \texttt{x}) := \mathcal{NN}_f(\texttt{f}) \odot \mathcal{NN}_k(\texttt{k}) \cdot \mathcal{NN}_G(\texttt{x})^\top \in \mathbb{R}^{M \times N}.
\end{equation}
In the formulation \eqref{eq0608_1}, $\texttt{k} \in \mathbb{R}^{M \times n_p}$ denotes the parameter input with $n_p$ being the number of parameter features. The output of the corresponding network is given by $\mathcal{NN}_k(\texttt{k}) \in \mathbb{R}^{M \times n}$. The symbol `$\odot$' denotes element-wise (Hadamard) product between two tensors of shape $M \times n$. All other notations are consistent with the definitions given above. The subsequent Figure \ref{fig0608_1} illustrates the specific network architecture.

Regardless of whether network \eqref{eq0607_6} or network \eqref{eq0608_1} is considered, the training is entirely data-driven. Specifically, to train network \eqref{eq0607_6}, a set of source terms $f$ is randomly generated, and the corresponding PDE is solved using a high-fidelity numerical solver to obtain the solution $u$. This process yields the necessary training data. Owing to the linearity of the differential operator, more than $N$ distinct training samples must be collected to ensure the network is sufficiently expressive. For the training of network \eqref{eq0608_1}, the parameter is uniformly sampled from its domain. For each sampled value of the parameter, the aforementioned procedure of generating $f$ and computing $u$ is repeated. The loss function used is the standard mean squared error (MSE) loss.

\begin{remark} \label{rmk3}
\textbf{(1) Relation to DeepONet.}
Although the design of network \eqref{eq0607_6} is motivated by the viewpoint of numerical integration in \eqref{eq0607_4} or \eqref{eq0607_5}, the resulting architecture turns out to be similar to that of DeepONet. However, once we incorporate the PDE-parameter input into the network, the structure \eqref{eq0608_1} becomes slightly different from the standard DeepONet.

\medskip

\textbf{(2) Generation of random source terms.}
In practice, the random source terms $f$ used for training are generated using two strategies:
\begin{itemize}
    \item \textbf{Gaussian-filtered noise:} Starting from standard Gaussian noise, apply Gaussian filters with different parameters to introduce smoothness.
    \item \textbf{Functional forms:} Construct source terms using random combinations of trigonometric basis functions and other functions. For instance,
    \[
    f(x,y) = a \cdot \sin(b\pi x) \cdot \cos(c\pi y), \quad \text{or similar variants,}
    \]
    where the coefficients $a, b, c$ are sampled uniformly from a predefined range.
\end{itemize}
These strategies ensure both smoothness and diversity in the source terms, thus enhancing the generalization capability of the trained model.
\end{remark}

\subsubsection{Operator Learning for Boundary Conditions} \label{Operator Learning for Boundary Conditions}
In this section, we consider the boundary integral equation \eqref{eq0607_3}, with the goal of learning a neural network that can directly predict the density function $\varphi$ from the Dirichlet boundary data $g^D$ and the parameters of the differential operator $\mathcal{L}$. Once the density function $\varphi$ is predicted, the solution $u_2$ to the PDE \eqref{eq0606_5} can be recovered via the representation formula \eqref{eq0607_2}.

Regarding the network architecture, since the input includes two parts—the boundary condition $g^D$ and the parameter of the operator $\mathcal{L}$—this neural network is adopted:
\begin{equation} \label{eq0608_2}
\mathcal{NN}_2(\texttt{k}, \texttt{g}^\texttt{D}) = \mathcal{NN}_{\text{output}}\left[\mathcal{NN}_k(\texttt{k}) \odot \mathcal{NN}_g (\texttt{g}^\texttt{D})\right],
\end{equation}
where $\odot$ denotes the element-wise (Hadamard) product. In this formulation, suppose the input batch size is $M$. Then, $\texttt{k} \in \mathbb{R}^{M \times n_p}$ represents the operator parameters, with $n_p$ being the number of parameter features, and $\texttt{g}^{\texttt{D}} \in \mathbb{R}^{M \times N_{\text{bd}}}$ denotes the Dirichlet boundary condition evaluated at $N_{\text{bd}}$ sampled points on the boundary. The sub-networks $\mathcal{NN}_k$ and $\mathcal{NN}_g$ map their respective inputs to feature representations in $\mathbb{R}^{M \times N_*}$, where the internal dimension $N_*$ need not equal $N_{\text{bd}}$, but must be consistent across the two branches to allow the Hadamard product. The output network $\mathcal{NN}_{\text{output}}$ then maps this product to a prediction for the density function $\varphi \in \mathbb{R}^{M \times N_{\text{bd}}}$. In this framework, the sub-networks $\mathcal{NN}_k$, $\mathcal{NN}_g$, and $\mathcal{NN}_{\text{output}}$ can be instantiated using architectures such as MLPs, ResNets, or other models.

For training, sample operator parameters $\texttt{k}$ and Dirichlet boundary functions $g^D$ from their respective spaces to form a dataset 
\[
\{(\texttt{k}_i, \texttt{g}^{\texttt{D}}_j) \mid i = 1,\dots,M_\lambda, \, j = 1,\dots,M_g \} = \{\texttt{k}_i\}_{i=1}^{M_\lambda} \times \{\texttt{g}^{\texttt{D}}_j\}_{j=1}^{M_g},
\]
where each $\texttt{k}_i \in \mathbb{R}^{n_p}$, $\texttt{g}^{\texttt{D}}_j \in \mathbb{R}^{N_{\text{bd}}}$ represents the values of the function $g^D_j$ at $N_{\text{bd}}$ sampled points $\{\mathbf{x}_l\}_{l=1}^{N_{\text{bd}}}$ on the boundary $\partial \Omega$, and `$\times$' means the Cartesian product. Let $\Theta$ denote the set of trainable parameters in the neural network $\mathcal{NN}_2$. The training objective is to minimize the discrepancy between the left-hand side of equation \eqref{eq0607_3} with predicted density function $\varphi^\Theta$ and the right-hand side of the boundary integral equation \eqref{eq0607_3}. 

Suppose that $\partial\Omega$ admits a 
$2\pi$-periodic parametrization $\gamma:[0,2\pi]\to\mathbb{R}^2$. 
For $\mathbf{x}=\gamma(s)$ and $\mathbf{y}=\gamma(t)$, the double-layer potential takes the form
\[
(\mathcal{D}\varphi)(\gamma(s)) 
= \int_{0}^{2\pi} K(s,t)\,\varphi(t)\,dt,
\]
where the kernel is
\begin{equation}\label{eq1119}
K(s,t)
= \frac{\partial G_0(\gamma(s),\gamma(t))}{\partial n_{y}}
\,|\gamma'(t)|.
\end{equation}
Although $\partial G_0(x,y)/\partial n_y$ is weakly singular as $x\to y$, 
the parametrization $\gamma$ ensures that the one-dimensional kernel $K(s,t)$ is 
nonsingular and admits a well-defined diagonal value $K(s,s)$; see~\cite{atkinson1997numerical}. 
Thus, a standard quadrature rule with nodes $\{t_l\}_{l=1}^{N_{\mathrm{bd}}}$ and 
weights $\{\omega_l\}_{l=1}^{N_{\mathrm{bd}}}$ yields the discretization
\[
(\mathcal{D}\varphi)(\gamma(s_j))
\approx 
\sum_{l=1}^{N_{\mathrm{bd}}} 
\omega_l\, K(s_j,t_l)\, \varphi(t_l).
\]

For the $i$-th training sample with PDE parameter $\texttt{k}_i$ and Dirichlet data 
$\texttt{g}^{\texttt{D}}_i$, the neural network predicts the density
\[
\varphi_i^\Theta
=
\mathcal{NN}_2(\texttt{k}_i, \texttt{g}^{\texttt{D}}_i;\Theta)\in\mathbb{R}^{1\times N_{\mathrm{bd}}}.
\]
The loss function enforcing the boundary integral equation is then defined as
\begin{equation}\label{eq:new_loss}
\text{Loss}_2(\Theta)
=
\frac{1}{M\,N_{\mathrm{bd}}}
\sum_{i=1}^{M}
\sum_{j=1}^{N_{\mathrm{bd}}}
\left|
\frac{1}{2}\,(\varphi_i^\Theta)_j
+
\sum_{l=1}^{N_{\mathrm{bd}}}
\omega_l\,K(s_j,t_l)\,(\varphi_i^\Theta)_l
-
(\texttt{g}^{\texttt{D}}_i)_j
\right|^2.
\end{equation}
This formulation incorporates the correct diagonal behavior of the double-layer kernel 
through the parametrized nonsingular kernel $K(s,t)$.

\begin{remark}
\textbf{(1) On the treatment of the singular kernel.}
The loss function~\eqref{eq:new_loss} is constructed from the boundary integral
equation involving the double-layer potential. 
Although the derivative of the fundamental solution 
$\partial G_0(x,y)/\partial n_y$ is singular in the limit $x\to y$, 
a smooth parametrization of the boundary $\partial\Omega$ transforms the boundary
integral into a one-dimensional integral with respect to the parameter~$t$.  
Under such a parametrization, the resulting kernel $K(s,t)$ is nonsingular 
and possesses a well-defined diagonal value, allowing standard quadrature 
without special singular treatment; see~\cite{atkinson1997numerical}.

\textbf{(2) Motivation for predicting only the density.}
In this work, only the boundary density $\varphi$ rather than 
the solution $u_2$ inside the domain is predicted.
This significantly reduces the output dimension in operator learning,
since both the Dirichlet data $g^D$ and the density $\varphi$ live on 
$\partial\Omega$. 
Furthermore, the training procedure is fully self-supervised, as the loss 
function is derived from the boundary integral equation itself, requiring
no solution labels.
Consequently, the data generation process only involves sampling random
boundary functions, resulting in both computational efficiency and improved 
generalization, as confirmed in numerical experiments. 
We further note that related work has also employed boundary integral theory to construct loss functions or neural network architectures, e.g., \citep{CSIAM-AM-4-275,lin2023bi}; however, the motivation, target problems, and technical implementation differ substantially from the present approach.
\end{remark}

\subsection{Solving Evolutionary PDEs by NEKM} \label{Solving Evolutionary PDEs by NEKM}
In this section, the heat equation, wave equation, and Schrödinger equation are taken as examples to illustrate in detail how the trained models can be utilized to solve time-evolution equations, thereby accomplishing the final realization of NEKM. It is worth noting that this method can also be extended to a broader class of equations, such as the Allen–Cahn equation, the Cahn–Hilliard equation. For instance, the Allen–Cahn equation can be decomposed into elliptic subproblems using techniques such as Strang \cite{li2022stability} or convex splitting schemes \cite{guan2014second}, while the Cahn–Hilliard equation can be tackled using boundary integral methods which are detailed in \cite{wei2020integral}.

\subsubsection{Heat Equations} \label{Heat Equations}
Consider the classical heat equation defined on a spatial domain $\Omega \subset \mathbb{R}^d$ with Dirichlet boundary conditions and some $T > 0$:
\begin{equation} \label{eq0608_4}
\begin{cases} 
u_t = \kappa \Delta u, & \textbf{x} \in \Omega, t \in [0, T], \\ 
u(\textbf{x},0) = u^0(\textbf{x}), & \textbf{x} \in \Omega, \\ 
u(\textbf{x},t) = g^D(\textbf{x}, t), & \textbf{x} \in \partial \Omega, t \in [0, T], 
\end{cases}
\end{equation}
where $\kappa > 0$ denotes the diffusion coefficient, $u^0$ is the initial condition, and $g^D$ is the time-dependent Dirichlet boundary condition.

To apply operator learning framework, we discretize the temporal derivative. For given $N \in \mathbb{N}_+$, consider a uniform time step $\tau = \frac{T}{N}$ and denote $u^n(\mathbf{x}) \approx u(\mathbf{x}, n\tau), g^n(\mathbf{x}) = g^D(\mathbf{x}, n\tau)$ for $n = 0, 1, \ldots, N$. The following two implicit time discretizations are considered:
\begin{itemize}
\item \textbf{Backward Euler scheme:}
\begin{equation} \label{eq0608_5}
\frac{u^{n+1} - u^n}{\tau} = \kappa \Delta u^{n+1},
\end{equation} 
which leads to the elliptic equation
\begin{equation} \label{eq0608_6}
\begin{cases} 
(\mathcal{I} - \kappa \tau \Delta) u^{n+1} = u^n, & \text{in}\ \Omega, \\ 
u^{n+1} = g^{n+1}, & \text{on}\ \partial\Omega, \\ 
\end{cases}
\end{equation}
where $\mathcal{I}$ means the identity operator.
\item \textbf{Crank–Nicolson scheme:}
\begin{equation} \label{eq0608_7}
\frac{u^{n+1} - u^n}{\tau} = \frac{1}{2} \kappa (\Delta u^{n} + \Delta u^{n+1}),
\end{equation}
which results in the elliptic equation
\begin{equation} \label{eq0608_8}
\begin{cases} 
(\mathcal{I} - \frac{\kappa \tau}{2} \Delta) u^{n+1} = (\mathcal{I} + \frac{\kappa \tau}{2} \Delta) u^n, & \text{in}\ \Omega, \\ 
u^{n+1} = g^{n+1}, & \text{on}\ \partial\Omega, \\ 
\end{cases}
\end{equation}
\end{itemize}
Under both schemes, the original time-dependent heat equation is reduced to solving a sequence of elliptic PDEs with Dirichlet boundary conditions at each time step. More importantly, each of these elliptic problems matches the structure of equation \eqref{eq0606_2}–\eqref{eq0606_3}, which means they can be directly solved using the trained operator models $\mathcal{NN}_1$ and $\mathcal{NN}_2$.
\begin{remark}
To avoid instability when computing the second-order spatial derivatives $u_{xx}$ and $u_{yy}$ in the right-hand side of the PDE in equation \eqref{eq0608_8}, we adopt the formulation $F^{j+1} \triangleq u^{j+1} + \frac{\kappa \tau}{2} \Delta u^{j+1} = 2u^{j+1} - F^j$, for $j = 0, 1, 2, \ldots$. In other words, $F^0$ can be computed exactly from the known initial condition $u^0$, while $F^1, F^2, \ldots$ can be obtained recursively using this relation. This recursive formulation for the right-hand side remains applicable to the subsequent numerical schemes and will not be discussed further.
\end{remark}

\subsubsection{Wave Equations} \label{Wave Equations}
Consider the classical wave equation with Dirichlet boundary conditions:
\begin{equation} \label{eq0608_wave_1}
\begin{cases} 
u_{tt} = \Delta u, & \textbf{x} \in \Omega \subset \mathbb{R}^d,\ t \in [0, T], \\ 
u(\textbf{x},0) = u_0(\textbf{x}), \quad u_t(\textbf{x}, 0) = u_1(\textbf{x}), & \textbf{x} \in \Omega, \\ 
u(\textbf{x},t) = g^D(\textbf{x}, t), & \textbf{x} \in \partial \Omega,\ t \in [0, T],
\end{cases}
\end{equation}
where $u^0$ and $u^1$ denote the initial displacement and velocity, respectively, and $g^D$ is the time-dependent Dirichlet boundary condition.

To apply operator learning framework, we discretize the temporal derivatives. For given $N \in \mathbb{N}_+$, consider a uniform time step $\tau = \frac{T}{N}$ and denote $u^n(\mathbf{x}) \approx u(\mathbf{x}, n\tau), g^n(\mathbf{x}) = g^D(\mathbf{x}, n\tau)$ for $n = 0, 1, \ldots, N$. The implicit $\theta$-scheme \cite{leveque2007finite} implies that:
\begin{equation} \label{eq0608_wave_2}
\frac{u^{n+1} - 2u^n + u^{n-1}}{\tau^2} = \theta \Delta u^{n+1} + (1 - 2\theta) \Delta u^n + \theta \Delta u^{n-1}, \quad \theta \in [0,1],
\end{equation}
which leads to the elliptic equation:
\begin{equation} \label{eq0608_wave_3}
\begin{cases}
(\mathcal{I} - \theta \tau^2 \Delta) u^{n+1} = 2u^n - u^{n-1} + \tau^2 \left[(1 - 2\theta)\Delta u^n + \theta \Delta u^{n-1} \right], & \text{in } \Omega, \\
u^{n+1} = g^{n+1}, & \text{on } \partial \Omega.
\end{cases}
\end{equation}
This formulation allows the wave equation to be solved as a sequence of elliptic equations with Dirichlet boundary conditions. This scheme becomes unconditionally stable and second-order accurate when $\theta \in [\frac{1}{4}, \frac{1}{2}]$. Similar to the heat equation, each step fits the form of equation \eqref{eq0606_2}–\eqref{eq0606_3} and can thus be solved using the trained operator models $\mathcal{NN}_1$ and $\mathcal{NN}_2$.

\subsubsection{The Schrödinger equation} \label{The Schrödinger equation}
The Schrödinger equation on a bounded irregular domain \(\Omega \subset \mathbb{R}^2\) with Dirichlet boundary conditions is given by:
\begin{equation} \label{eq0608_S1}
    -i\partial_t u(\mathbf{x},t) + \Delta u(\mathbf{x},t) + v(\mathbf{x})u(\mathbf{x},t) + w|u(\mathbf{x},t)|^2 u(\mathbf{x},t) = 0 \quad \text{in}\ \Omega \times (0,T],
\end{equation}
where \(v(\mathbf{x}) \in \mathbb{R}\) is a prescribed potential function, \(w \geq 0\) is a nonlinear coefficient, and \(u(\mathbf{x},t) \in \mathbb{C}\) is the unknown function. The initial and Dirichlet boundary conditions are:
\begin{equation} \label{eq0608_S2_1}
    u(\mathbf{x},0) = u_0(\mathbf{x}), \quad \mathbf{x} \in \Omega,
\end{equation}
\begin{equation} \label{eq0608_S2_2}
    u(\mathbf{x},t)\big|_{\partial \Omega} = g^D(\mathbf{x},t), \quad t > 0.
\end{equation}

To apply operator learning framework, we discretize the time derivative. For given $N \in \mathbb{N}_+$, consider a uniform time step $\tau = \frac{T}{N}$ and denote $u^n(\mathbf{x}) \approx u(\mathbf{x}, n\tau), g^n(\mathbf{x}) = g^D(\mathbf{x}, n\tau)$ for $n = 0, 1, \ldots, N$. The Schrödinger equation can be discretized using the Strang splitting method or the Lie-Trotter splitting method. The equation can be split into linear and nonlinear subproblems:

\noindent 1. \textbf{Linear subproblem:}
\begin{equation} \label{eq0608_S3}
    i\partial_t u = \Delta u,
\end{equation}
with solution operator \(\mathcal{S}^A\).

\noindent 2. \textbf{Nonlinear subproblem:}
\begin{equation} \label{eq0608_S4}
    i\partial_t u = vu + w|u|^2u,
\end{equation}
with solution operator \(\mathcal{S}^B\).

Consider the following two implicit time discretizations:
\begin{itemize}
\item \textbf{Strang splitting:}
The Strang splitting scheme combines these operators:
\begin{equation} \label{eq0608_S5}
    u^{n+1} = \mathcal{S}^A\left(\frac{\tau}{2}\right) \mathcal{S}^B(\tau) \mathcal{S}^A\left(\frac{\tau}{2}\right) u^n.
\end{equation}
Introducing two intermediate variables, $u^{*}(x)$ and $u^{**}(x)$, the Strang splitting method can be employed to solving the problem \eqref{eq0608_S1}, which leads to the following algorithm:
\begin{equation} \label{eq0608_S6}
\left\{\begin{array}{l}
u^*(\mathbf{x})=u^n(\mathbf{x})-i \frac{\tau}{2} \Delta u^n(\mathbf{x}) \\
u^{* *}(\mathbf{x})+\frac{i \tau}{2}\left(v(\mathbf{x})+w\left|u^{* *}(\mathbf{x})\right|^2\right) u^{* *}(\mathbf{x})=u^*(\mathbf{x})-\frac{i \tau}{2}\left(v(\mathbf{x})+w\left|u^*(\mathbf{x})\right|^2\right) u^*(\mathbf{x}) \\
u^{n+1}(\mathbf{x})+i \frac{\tau}{2} \Delta u^{n+1}(\mathbf{x})=u^{* *}(\mathbf{x})
\end{array}\right.
\end{equation}
Note that the first and third equations of \eqref{eq0608_S6} are obtained by employing the forward Euler scheme and backward Euler scheme to \eqref{eq0608_S3}, respectively. Moreover, the second equation of \eqref{eq0608_S6} is obtained by employing the Crank-Nicolson scheme to \eqref{eq0608_S4}.
 
\item \textbf{Lie-Trotter splitting:}
Moreover, Lie-Trotter splitting suggests the following procedure:
\begin{equation} \label{eq0608_S7}
    u^{n+1} = \mathcal{S}^A(\tau) \mathcal{S}^B\left(\tau\right) u^n.
\end{equation}
Introducing one intermediate variable $u^{*}(x)$, the Lie-Trotter splitting method can be employed to solve the problem \eqref{eq0608_S1}, which leads to the following algorithm:
\begin{equation} \label{eq0608_S8}
\left\{\begin{array}{l}
u^{* }(\mathbf{x}) + \frac{i \tau}{2} \left(v(\mathbf{x})+w\left|u^{*}(\mathbf{x})\right|^2\right) u^{*}(\mathbf{x})=u^n(\mathbf{x})-\frac{i \tau}{2} \left(v(\mathbf{x})+w\left|u^n(\mathbf{x})\right|^2\right) u^n(\mathbf{x}) \\
u^{n+1}(\mathbf{x})+i \tau \Delta u^{n+1}(\mathbf{x})=u^{*}(\mathbf{x})
\end{array}\right.
\end{equation}
\end{itemize}
The second equation in \eqref{eq0608_S6} and the first equation in \eqref{eq0608_S8} are essentially nonlinear equations defined at each spatial point $\mathbf{x} \in \Omega$, which can be solved using Newton's method. The first equation in \eqref{eq0608_S6} can be computed directly. The remaining equations—the third equation in \eqref{eq0608_S6} and the second equation in \eqref{eq0608_S8}—together with the associated Dirichlet boundary conditions, constitute a generalization of equations \eqref{eq0606_2}–\eqref{eq0606_3} in the case where $u \in \mathbb{C}$—or equivalently, a PDE system for $\text{Re}(u)$ and $\text{Im}(u)$. For this type of PDE system, a similar approach as in Section \ref{Operator Learning} can be followed: first derive the corresponding boundary integral equations, and then perform operator learning for the source term $\vec{f}$ and the Dirichlet boundary condition $\vec{g^D}$. The implementation details are provided in Appendix \ref{Operator Learning for PDE Systems}. The Schrödinger equation \eqref{eq0608_S1}–\eqref{eq0608_S2_2} can then be solved using the trained operator models.

\section{Numerical Experiments} \label{Numerical Experiments}
This section presents numerical results for several canonical equations, including the modified Helmholtz equation, the heat equation, the wave equation, and the Schr\"{o}dinger equation. Specifically, operator learning is first demonstrated for the modified Helmholtz equation on a rectangular domain, and the trained model is subsequently applied to solve the heat and wave equations. Next, following the same approach, operator learning is performed for the discretized PDE system derived from the Schr\"{o}dinger equation, and the learned models are used to compute its solution. To further demonstrate the generality of the method, a more complex computational domain—a petal-shaped domain—is considered; on this domain, the model is retrained and applied to solve the heat equation. All experiments yield results with high accuracy, indicating strong generalization capability and precision of the approach. In addition, an example applying the model to uncertainty quantification is included, showcasing the extensibility of the method. While some details of the network architecture and training procedures are not fully elaborated in the text, the complete implementation code for all experiments is available at \url{https://github.com/vstppen/NEKM2}.

\subsection{Modified Helmholtz Equations} \label{Modified Helmholtz Equations}
In this section, the domain is $\Omega = [0, 1]^2$ and the modified Helmholtz equation is given by:
\begin{equation} \label{eq1127_3} 
\begin{cases}
\Delta u - \frac{1}{\kappa} u = 0, & \mathbf{x} \in \Omega, \\
u = g, & \mathbf{x} \in \partial \Omega,
\end{cases}
\end{equation}
where $\kappa > 0$ is a given parameter. The operator $\mathcal{L} = \Delta - \frac{1}{\kappa}\mathcal{I}$ (with $\mathcal{I}$ the identity operator) is uniformly elliptic, coercive and positive definite. We adopt the operator learning framework described in Section~\ref{Operator Learning for Boundary Conditions} and train two separate models 
$\mathcal{NN}_2(\texttt{k}, \texttt{g}^\texttt{D})$: \textbf{Model 1} for $\kappa \in [0.05, 0.1]$ and \textbf{Model 2} for $\kappa \in [\frac{1}{128}, \frac{1}{32}]$. This model partition is introduced to balance the approximation accuracy across different parameter regimes, since a single model trained over the entire range of $\kappa$ may suffer from degraded performance or exhibit significantly lower accuracy in certain subranges. The training dataset consists of $M_g = 10{,}000$ samples of boundary data $g$ and $M_{\lambda} = 11$ samples of $\kappa$, with boundary values discretized over $N_{\textbf{bd}} = 512$ points. In numerical experiments, we adopt the neural network \eqref{eq0608_2} implemented as follows. 
The sub-network $\mathcal{NN}_k$ takes the scalar parameter $\tau$ as input and is realized as an MLP with two hidden layers of width $3 N_{\mathrm{bd}}/4$ and ReLU activation, producing a feature vector in $\mathbb{R}^{3N_{\mathrm{bd}}/4}$. 
The sub-network $\mathcal{NN}_g$ processes the boundary data $g \in \mathbb{R}^{N_{\mathrm{bd}}}$ through a single linear layer, yielding a feature vector in $\mathbb{R}^{3N_{\mathrm{bd}}/4}$. Given the outputs of the two subnetworks, denoted as the $\tau$-features and the $g$-features, the final output is obtained by taking their elementwise product and applying a final linear transformation, resulting in a vector in $\mathbb{R}^{N_{\mathrm{bd}}}$. 

\begin{figure}[htbp]
    \centering
    \includegraphics[width=0.82\textwidth]{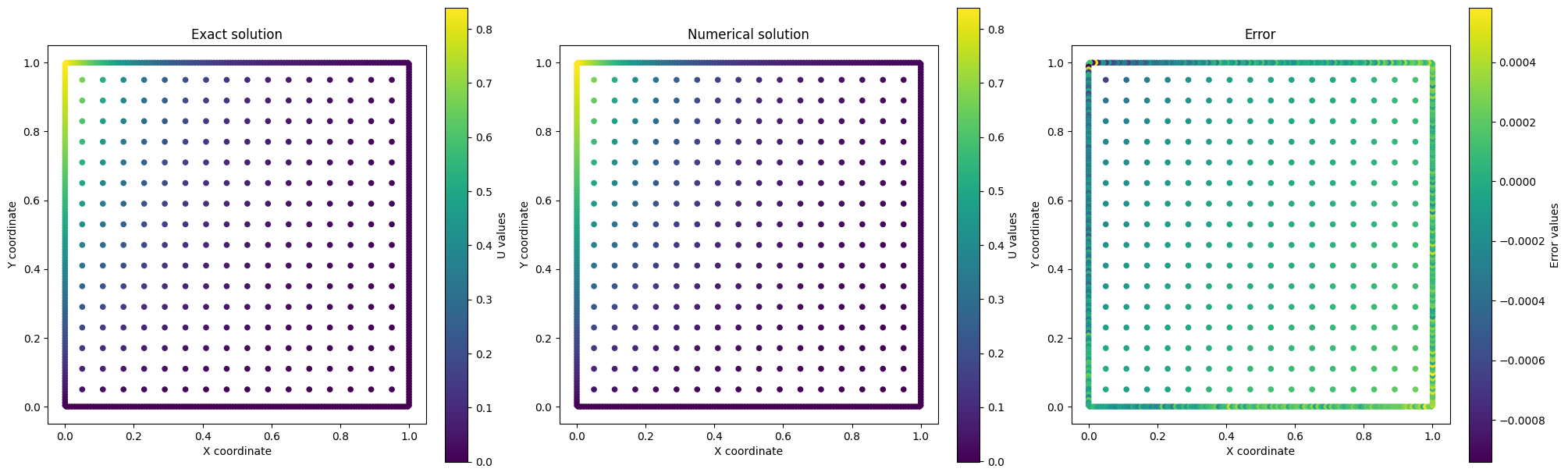} 
    \caption{Numerical results for Model 1 with $\kappa = 0.067$ which is not in the training set.}
    \label{fig1127_1}
\end{figure}

\begin{figure}[htbp]
    \centering
    \includegraphics[width=0.82\textwidth]{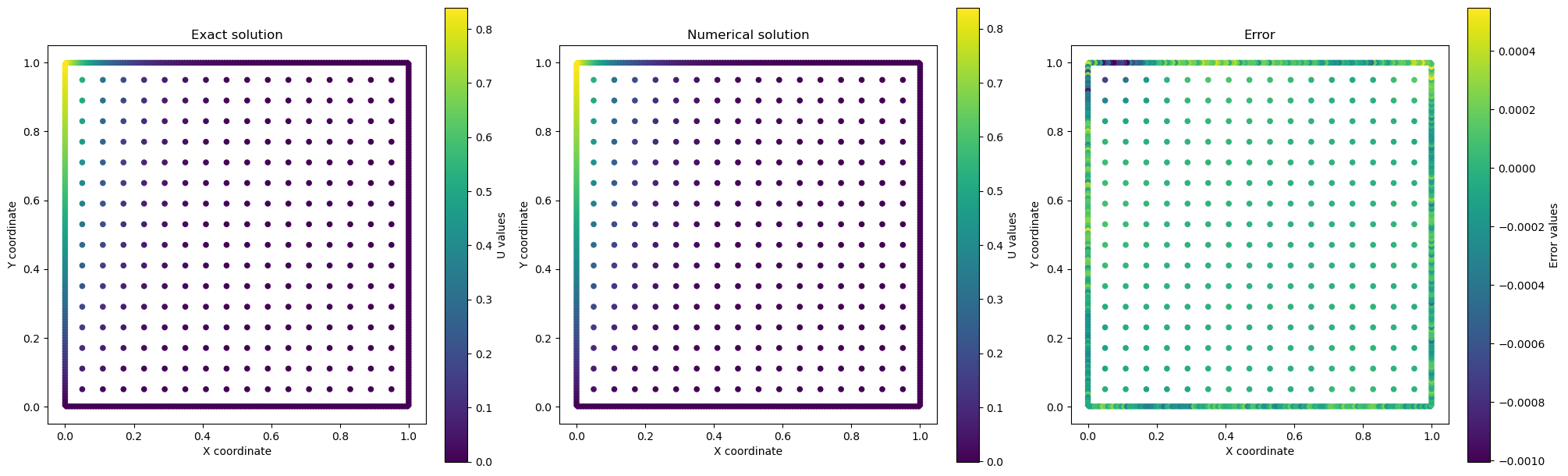} 
    \caption{Numerical results for Model 2 with $\kappa = \frac{1}{72}$ which is not in the training set.}
    \label{fig0219}
\end{figure}

To evaluate the accuracy of the trained models, we use exact solutions of the form $$u_{\kappa}(x, y) = e^{-\sqrt{1 + \frac{1}{\kappa}} x} \sin(y).$$ The evaluation is performed on a uniform \(16 \times 16\) Cartesian grid in the interior region \([0.05,\,0.95]^2 \subset \bar{\Omega}\). Selected numerical results for \textbf{Model 1} are presented in Table~\ref{tab1127_1} and Figure \ref{fig1127_1}, while those for \textbf{Model 2} are provided in Table~\ref{tab0219} and Figure \ref{fig0219}.
\begin{table}[htbp]
    \centering
    \begin{tabular}{|c|c|c|c|c|} \hline 
         $\kappa$ &  Absolute $L_2$ error&  Absolute $L_{\infty}$ error&  Relative $L_2$  error& Relative $L_{\infty}$ error\\ \hline 
         0.05& 1.777E-4  & 6.045E-4 & 0.1120\% & 0.09340\% \\ \hline 
         0.067& 1.112E-4 & 6.457E-4 & 0.06457\% & 0.09691\% \\ \hline 
         0.1& 1.495E-4 & 5.506E-4 & 0.07803\% & 0.07989\% \\ \hline 
    \end{tabular}
    \caption{Numerical errors for Model 1.}
    \label{tab1127_1}
\end{table}

\begin{table}[htbp]
    \centering
    \begin{tabular}{|c|c|c|c|c|} \hline 
         $\kappa$ &  Absolute $L_2$ error&  Absolute $L_{\infty}$ error&  Relative $L_2$  error& Relative $L_{\infty}$ error\\ \hline 
         $\frac{1}{128}$& 1.064E-4  & 7.870E-4 & 0.1249\% & 0.1707\% \\ \hline 
         $\frac{1}{72}$ & 6.563E-5 & 6.880E-4 & 0.06209\% & 0.1296\% \\ \hline 
         $\frac{1}{32}$ & 1.171E-4 & 3.783E-4 & 0.08495\% & 0.06199\% \\ \hline 
    \end{tabular}
    \caption{Numerical errors for Model 2.}
    \label{tab0219}
\end{table}

Next, we evaluate our operator learning approach on the modified Helmholtz equation with nonzero source term:
\begin{equation} \label{eq0618_1}
\begin{cases}
\Delta u - \frac{1}{\kappa} u = f, & \mathbf{x} \in \Omega, \\
u = 0, & \mathbf{x} \in \partial \Omega,
\end{cases}
\end{equation}
To generate the training data, a variety of source terms \( f \) are constructed. One part of the dataset was generated by sampling random fields using \texttt{np.random.randn}, followed by Gaussian smoothing to control regularity. The other part consisted of analytical expressions formed by linear combinations of trigonometric functions such as 
\[
\sin(k_1 x) \sin(k_2 y), \quad \sin(k_3 x) \cos(k_4 y), \quad \cos(k_5 x) \sin(k_6 y), \quad \cos(k_7 x) \cos(k_8 y),
\]
where \( k_1, \dots, k_8 \) are randomly selected coefficients. For each sampled $\kappa$, a sufficient number of source terms was generated, and the corresponding solutions \( u \) were computed using classical numerical solvers. Based on the framework introduced in Section~\ref{Operator Learning for Source Terms}, we trained two models $\mathcal{NN}_1(\texttt{f}, \texttt{k}, \texttt{x})$: \textbf{Model 3} for $\kappa \in [0.05, 0.1]$ and \textbf{Model 4} for $\kappa \in [\frac{1}{128}, \frac{1}{32}]$, where $\kappa$ values were uniformly sampled within the respective intervals. In our numerical experiments, $N_{\text{bd}} = 41$ is the grid size of a uniform Cartesian grid on $\Omega$, corresponding to the parameters setting $N = n = N_{\text{bd}}^2$ in Section \ref{Operator Learning for Source Terms}. The operator $\mathcal{NN}_1$ in \eqref{eq0608_1} is implemented by the following. The sub-network $\mathcal{NN}_k$ takes the parameter $\tau$ as input and is realized as a four-layer MLP, each hidden layer having width $N$, with ReLU activation, and resulting a feature vector of $\mathbb{R}^n$. The sub-network $\mathcal{NN}_G$ maps all grid points $\mathbf{x} \in \mathbb{R}^{N \times 2}$ to a feature tensor in $\mathbb{R}^{N \times n}$ through an MLP of identical depth and width in $\mathcal{NN}_k$. The sub-network $\mathcal{NN}_f$ is simply the identity mapping since $N = n$. 

To assess the performance of the models, we employed the exact solution 
\[
u(x, y) = x(1 - x) y(1 - y) \exp(0.6x + 0.8y)
\]
and computed the associated source term \( f \) for different values of $\kappa$ to serve as test inputs. The evaluation is conducted on the uniform $41 \times 41$ Cartesian grid in $\bar{\Omega}$. Unless otherwise specified, all subsequent error evaluations will use this same set of test points. The numerical results for \textbf{Model 3} are reported in Table~\ref{tab1127_2}, and those for \textbf{Model 4} are provided in Table~\ref{tab0219_2}.
\begin{table}[htbp]
\centering
\begin{tabular}{cccccc}
\toprule
\textbf{$\kappa$} & \textbf{Abs. L2 Norm} & \textbf{Abs. Max Norm} & \textbf{Rel. L2 Norm} & \textbf{Rel. Max Norm} \\
\midrule
0.055 &    7.139E-5 &     2.571E-4  &   0.1052\%     &     0.1921\% \\
0.065   &  6.736E-5  &    2.298E-4   &  0.09926\%     &     0.1717\% \\
0.075   &  6.486E-5  &    2.189E-4   &   0.09559\%     &     0.1636\% \\
0.085 &    6.602E-5     & 2.436E-4  &    0.09728\%    &      0.1820\% \\
0.095  &   6.893E-5   &   2.547E-4   &  0.1016\%    &       0.1903\% \\
\bottomrule
\end{tabular}
\caption{Numerical errors for \textbf{Model 3}.}
\label{tab1127_2}
\end{table}

\begin{table}[htbp]
\centering
\begin{tabular}{cccccc}
\toprule
\textbf{$\kappa$} & \textbf{Abs. L2 Norm} & \textbf{Abs. Max Norm} & \textbf{Rel. L2 Norm} & \textbf{Rel. Max Norm} \\
\midrule
0.009 & 2.532E-4 & 1.436E-4 & 0.05827\% & 0.1073\% \\
0.013 & 2.329E-4 & 1.245E-4 & 0.05360\% & 0.1065\% \\
0.017 & 2.313E-4 & 1.501E-4 & 0.05322\% & 0.1121\% \\
0.021 & 2.465E-4 & 1.736E-4 & 0.05674\% & 0.1297\% \\
0.025 & 2.484E-4 & 1.654E-4 & 0.05716\% & 0.1358\%  \\
0.029 & 2.406E-4 & 1.549E-4 & 0.05538\% & 0.1158\% \\
\bottomrule
\end{tabular}
\caption{Numerical errors for \textbf{Model 4}.}
\label{tab0219_2}
\end{table}

The numerical results in the tables indicate that the method achieves uniformly low errors across a wide range of parameter values, including those beyond the training range, thus demonstrating its strong interpolation and extrapolation capabilities. In particular, both the $L^2$ and $L^\infty$ relative errors remain below 0.2\% within the training interval, and the errors remain acceptably small even for extrapolated values of $\kappa$. This highlights the high approximation accuracy of the learned operator models. Such high-fidelity performance is a key reason why the trained networks can be further employed as evolutionary kernels in time-dependent PDEs. Since time discretization schemes often lead to sequences of elliptic problems with varying source terms and boundary conditions, the ability of NEKM to produce accurate and stable solutions across these inputs makes it particularly suitable for time evolution tasks, as will be demonstrated in the following sections.

\medskip
\noindent \textbf{Comparison with DeepONet.} Now, this following modified Helmholtz equation is considered:
\begin{equation} \label{eq0721_1} 
\begin{cases}
\Delta u - \frac{1}{\kappa} u = f, & \mathbf{x} \in \Omega, \\
u = g, & \mathbf{x} \in \partial \Omega,
\end{cases}
\end{equation}
To demonstrate the high accuracy of our method, we employ supervised learning to train three DeepONet-based networks, denoted as \textbf{Model D1}, \textbf{Model D3}, and \textbf{Model D}.
\begin{itemize}
\item \textbf{Model D1} has the similar functionality as \textbf{Model 1} and is used to solve equation \eqref{eq1127_3}. The input to the branch network in this DeepONet includes the equation parameter $\kappa \in [0.05, 0.1]$, the boundary condition $g$, while the trunk network takes the spatial coordinate $x$ as input. In our experiments, both the branch and trunk networks are four-layer MLPs with width $N^2 // 4$ (where $N = 41$, and $//$ denotes integer division), and output dimension also $N^2 // 4$, using ReLU activation. 

\item \textbf{Model D3} has the same functionality as \textbf{Model 3} and is used to solve equation \eqref{eq0618_1}. The input to the branch network in this DeepONet includes the equation parameter $\kappa \in [0.05, 0.1]$, the source term $f$, while the trunk network takes the spatial coordinate $x$ as input. In our experiments, both the branch and trunk networks are four-layer MLPs with width $N^2 // 4$ (where $N = 41$, and $//$ denotes integer division), and output dimension also $N^2 // 4$, using ReLU activation. 

\item \textbf{Model D} is designed to solve equation \eqref{eq0721_1}. The input to the branch network in this DeepONet includes the equation parameter $\kappa \in [0.05, 0.1]$, the source term $f$, and the boundary condition $g$, while the trunk network takes the spatial coordinate $x$ as input. In our experiments, both the branch and trunk networks are four-layer MLPs with width $3 N^2 // 4$ (where $N = 41$, and $//$ denotes integer division), and output dimension also $3 N^2 // 4$, using ReLU activation. Such a choice of network parameters ensures that the resulting DeepONet has approximately the same number of parameters as the total of \textbf{Model 1} and \textbf{Model 3} previously trained.
\end{itemize}

To generate training data $f$ and $g$, we adopt the schemes discussed in \cref{rmk3}, and additionally employ the Gaussian random field method introduced in the DeepONet paper \cite{lu2021learning} to generate function values on the grid points of $[0,1]$. For boundary data, only values on $\partial \Omega$ are used. The comparison between \textbf{Model 1} and \textbf{Model D1} is shown in Table \ref{tab1121_1}. Note that the selected test ground-truth solution is $u_{\kappa}(x, y) = e^{-\sqrt{1 + \frac{1}{\kappa}} x} \sin(y)$.
\begin{table}[htbp]
    \centering
    \setlength{\tabcolsep}{10pt} 
    \renewcommand{\arraystretch}{1.25} 
    \begin{tabular}{|c|l|c|c|c|} \hline
        \textbf{$\kappa$} & \textbf{Error type} & \textbf{Model 1} & \textbf{Model D1} & \textbf{Reduction (\%)} \\ \hline
        \multirow{2}{*}{0.05} 
        & Abs $L_2$ error        & 1.777E-4 & 3.833E-3 & 95.36 \\ 
        & Abs $L_\infty$ error  & 6.045E-4 & 1.331E-2 & 95.45 \\
        \hline
        \multirow{2}{*}{0.067} 
        & Abs $L_2$ error        & 1.112E-4 & 2.017E-3 & 94.50 \\ 
        & Abs $L_\infty$ error  & 6.457E-4 & 8.249E-3 & 92.20 \\
        \hline
        \multirow{2}{*}{0.1} 
        & Abs $L_2$ error        & 1.495E-4 & 3.818E-3 & 96.08 \\ 
        & Abs $L_\infty$ error  & 5.506E-4 & 2.315E-2 & 97.62 \\
        \hline
        \multicolumn{2}{|l|}{sample size for $g$} & 10000 & 25000 & -- \\ \hline
        \multicolumn{2}{|l|}{sample size for $u$} & 0 & 25000 & -- \\ \hline
        \multicolumn{2}{|l|}{training epochs} & 50000 & 100000 & -- \\ \hline
    \end{tabular}
    \caption{Comparison for \textbf{Model 1} and \textbf{Model D1}. The reduction is computed relative to Model D1.}
    \label{tab1121_1}
\end{table}

The comparison between \textbf{Model 3} and \textbf{Model D3} is shown in Table \ref{tab1121_2}. Note that the selected test ground-truth solution is $u(x, y) = x(1 - x) y(1 - y) \exp(0.6x + 0.8y)$. Four models in Tables \ref{tab1121_1} and \ref{tab1121_2} have been fully trained.
\begin{table}[htbp]
    \centering
    \setlength{\tabcolsep}{10pt} 
    \renewcommand{\arraystretch}{1.25} 
    \begin{tabular}{|c|l|c|c|c|} \hline
        \textbf{$\kappa$} & \textbf{Error type} & \textbf{Model 3} & \textbf{Model D3} & \textbf{Reduction (\%)} \\ \hline
        \multirow{2}{*}{0.055} 
        & Abs $L_2$ error        & 7.139E-5 & 1.765E-3 & 95.96 \\ 
        & Abs $L_\infty$ error  & 2.571E-4 & 4.902E-3 & 94.76 \\
        \hline
        \multirow{2}{*}{0.075} 
        & Abs $L_2$ error        & 6.486E-5 & 4.022E-3 & 98.39 \\ 
        & Abs $L_\infty$ error  & 2.189E-4 & 7.641E-3 & 97.14 \\
        \hline
        \multirow{2}{*}{0.095} 
        & Abs $L_2$ error        & 6.893E-5 & 3.148E-3 & 97.81 \\ 
        & Abs $L_\infty$ error  & 2.547E-4 & 7.538E-3 & 96.62 \\
        \hline
        \multicolumn{2}{|l|}{sample size for $f$} & 22000 & 100000 & -- \\ \hline
        \multicolumn{2}{|l|}{sample size for $u$} & 22000 & 100000 & -- \\ \hline
        \multicolumn{2}{|l|}{training epochs} & 50000 & 50000 & -- \\ \hline
    \end{tabular}
    \caption{Comparison for \textbf{Model 3} and \textbf{Model D3}. The reduction is computed relative to Model D3.}
    \label{tab1121_2}
\end{table}

Tables~\ref{tab1121_1} and~\ref{tab1121_2} demonstrate that the operator learning networks specifically designed in this paper for solving Equations~\eqref{eq0606_4} and~\eqref{eq0606_5} achieve significantly higher accuracy compared to the conventional DeepONet, while maintaining comparable training cost. An often overlooked issue, however, is whether the overall error in the combined solution of Equations~\eqref{eq0606_2}--\eqref{eq0606_3}, obtained by summing the separately computed solutions of Equations~\eqref{eq0606_4} and~\eqref{eq0606_5}, is dominated by one component or whether both components contribute similarly. To investigate this, the following experiments were conducted: solutions were computed using \textbf{Model D1} \& \textbf{Model D3}, \textbf{Model D1} \& \textbf{Model 3}, \textbf{Model 1} \& \textbf{Model D3}, and \textbf{Model 1} \& \textbf{Model 3}. By comparing errors from the last three methods with the benchmark obtained by solving Equations~\eqref{eq0606_4} and~\eqref{eq0606_5} separately using DeepONet, it is possible to identify whether the improvements in solving Equation~\eqref{eq0606_4} or Equation~\eqref{eq0606_5} contribute more significantly to the overall accuracy. For clarity, \textbf{Model 1} and \textbf{Model D1} are used for Equation~\eqref{eq0606_5}, while \textbf{Model 3} and \textbf{Model D3} are used for Equation~\eqref{eq0606_4}. Comparisons of the accuracy of these approaches are shown in Table \ref{tab1121_3} and Table \ref{tab0129}. Additionally, the results obtained by directly applying DeepONet to Equations~\eqref{eq0606_2}--\eqref{eq0606_3} (\textbf{Model D}) are also presented in the tables.
\begin{table}[htbp]
    \centering
    \resizebox{\textwidth}{!}{
    \setlength{\tabcolsep}{10pt} 
    \renewcommand{\arraystretch}{1.25} 
    \begin{tabular}{|c|c|c|c|c|c|} \hline
        \textbf{$\kappa$} & \textbf{Model D} & \textbf{Model D1 \& D3} & \textbf{Model D1 \& 3} & \textbf{Model 1 \& D3} & \textbf{Model 1 \& 3} \\ \hline
        0.05 & 0.069359 & 0.005380 & 0.003164 & 0.003258 & 0.000313 \\ \hline
        0.06 & 0.068182 & 0.004333 & 0.003119 & 0.002205 & 0.000201 \\ \hline
        0.07 & 0.066972 & 0.008474 & 0.003918 & 0.004951 & 0.000210 \\ \hline
        0.08 & 0.065781 & 0.007978 & 0.004156 & 0.004273 & 0.000179 \\ \hline
        0.09 & 0.064615 & 0.006582 & 0.003702 & 0.003740 & 0.000144 \\ \hline
        0.10 & 0.063473 & 0.005664 & 0.003538 & 0.004063 & 0.000216 \\  \hline
    \end{tabular}}
    \caption{Comparison of absolute $L_2$ errors with different $\kappa$ values using different models. The source term and Dirichlet boundary condition are given by modified Helmholtz equation \eqref{eq0721_1} with the exact solution $u(x, y) = \sin{(\frac{1}{\sqrt{2}}x)} \cos{(\frac{1}{\sqrt{2}}x)}$.}
    \label{tab1121_3}
\end{table}

\begin{table}[htbp]
    \centering
    \resizebox{\textwidth}{!}{
    \setlength{\tabcolsep}{10pt} 
    \renewcommand{\arraystretch}{1.25} 
    \begin{tabular}{|c|c|c|c|c|c|} \hline
        \textbf{$\kappa$} & \textbf{Model D} & \textbf{Model D1 \& D3} & \textbf{Model D1 \& 3} & \textbf{Model 1 \& D3} & \textbf{Model 1 \& 3} \\ \hline
        0.05 & 0.047516 & 0.005124 & 0.002957 & 0.004095 & 0.000222 \\ \hline
        0.06 & 0.046682 & 0.005869 & 0.002867 & 0.003745 & 0.000141 \\ \hline
        0.07 & 0.045906 & 0.006572 & 0.003188 & 0.003734 & 0.000148 \\ \hline
        0.08 & 0.045201 & 0.006146 & 0.003209 & 0.003490 & 0.000125 \\ \hline
        0.09 & 0.044508 & 0.005474 & 0.002954 & 0.003588 & 0.000100 \\ \hline
        0.10 & 0.043822 & 0.004745 & 0.003624 & 0.003601 & 0.000154 \\  \hline
    \end{tabular}}
    \caption{Comparison of absolute $L_2$ errors with different $\kappa$ values using different models. The source term and Dirichlet boundary condition are given by modified Helmholtz equation \eqref{eq0721_1} with the exact solution $u(x, y) = \sin{(ax)} \cos{(bx)}$ for $a = 0.5$ and $b = \sqrt{1 - a^2}$.}
    \label{tab0129}
\end{table}

From the results in Tables~\ref{tab1121_3} and Table~\ref{tab0129}, it can be observed that replacing only one of the components in \textbf{Model 1} \& \textbf{Model 3} with the operator learning network developed in this paper results in a comparable reduction in error. This indicates that, when using the DeepONet method to solve Equations~\eqref{eq0606_4} and~\eqref{eq0606_5} separately, neither component's error dominates the other. Furthermore, the results in Tables~\ref{tab1121_1} and~\ref{tab1121_2} show that the improvement of \textbf{Model 1} over \textbf{Model D1} is similar to that of \textbf{Model 3} over \textbf{Model D3}. Therefore, partitioning Equations~\eqref{eq0606_2}--\eqref{eq0606_3} into two subequations and training separate operator learning networks for each by the strategies in Section \ref{Operator Learning for Source Terms} and \ref{Operator Learning for Boundary Conditions} is a reasonable and effective strategy.

\subsection{Heat Equations} \label{Numerical Heat Equations} 
In this section, consider the following heat equation, where \( a, b \in (0, 1) \) are two constants satisfying \( a^2 + b^2 = 1 \):
\begin{equation} \label{eq6}
\begin{cases} 
u_t = \Delta u, & \mathbf{x} \in \Omega = [0, 1]^2,\ t \in [0, T], \\ 
u(\mathbf{x}, 0) = \sin(a x_1) \cos(b x_2), & \mathbf{x} = (x_1, x_2) \in \Omega, \\ 
u(\mathbf{x}, t) = e^{-t} \sin(a x_1) \cos(b x_2), & \mathbf{x} \in \partial \Omega,\ t \in [0, T].
\end{cases}
\end{equation}
The exact solution is given by
\[
u(\mathbf{x}, t) = e^{-t} \sin(a x_1) \cos(b x_2),
\]
which satisfies both the initial and boundary conditions. We evaluate the performance of the NEKM method in terms of accuracy for individual PDEs, computational efficiency for solving multiple PDEs in parallel, and effectiveness in uncertainty quantification.

\begin{itemize}
\item \textbf{Solving a single PDE.} In this part, we consider the case where \( a = b = \frac{1}{\sqrt{2}} \). The PDE is solved using the backward Euler scheme described in Section \ref{Heat Equations}, with the pre-trained \textbf{Model 1} and \textbf{Model 3} from Section \ref{Modified Helmholtz Equations} used to predict the solution at each time step. Representative results are presented in Table \ref{tab1127_4_1} and Figure \ref{fig1127_4_1}.
\begin{table}[htbp]
    \centering
    \begin{tabular}{|c|c|c|c|c|} \hline 
         $\tau$ &  Absolute $L_2$ error&  Absolute $L_{\infty}$ error&  Relative $L_2$  error& Relative $L_{\infty}$ error\\ \hline 
         0.1& 0.000558 & 0.001190 & 0.2314\% & 0.3281\% \\ \hline 
         $\frac{1}{15}$ & 0.000377 & 0.000939 & 0.1566\% & 0.2591\% \\ \hline 
         0.05 & 0.000271 & 0.000868 & 0.1123\% & 0.2394\% \\ \hline 
    \end{tabular}
    \caption{Numerical errors at time $T = 1$ for heat equation \eqref{eq6} with different time step $\tau$.}
    \label{tab1127_4_1}
\end{table}
\begin{figure}[htbp]
    \centering
    \begin{subfigure}
        \centering
        \includegraphics[width=\textwidth]{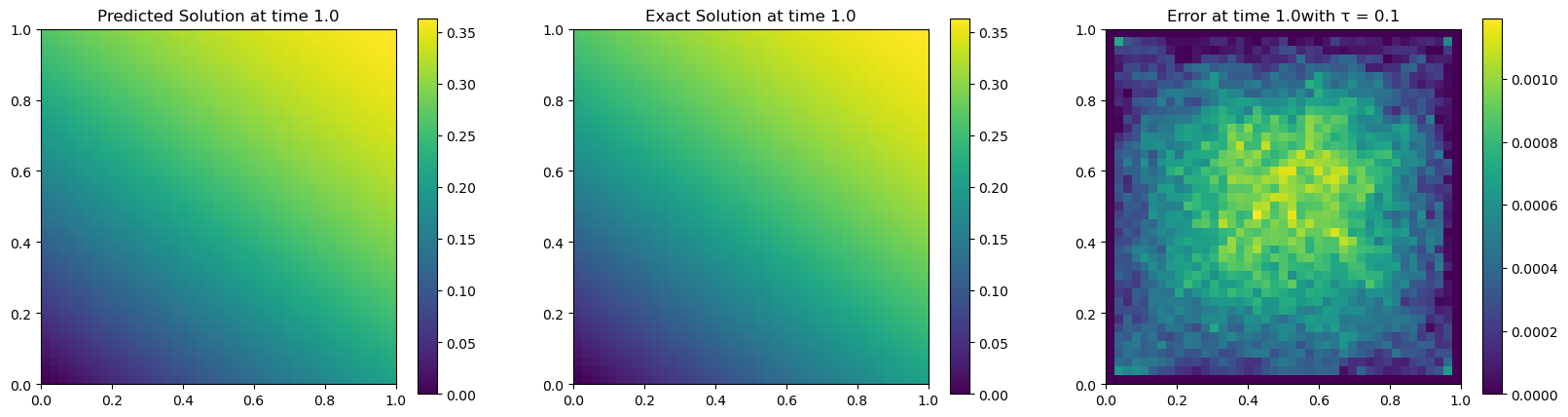}
        \label{fig1127_4_1_a}
    \end{subfigure}
    \begin{subfigure}
        \centering
        \includegraphics[width=0.95\textwidth]{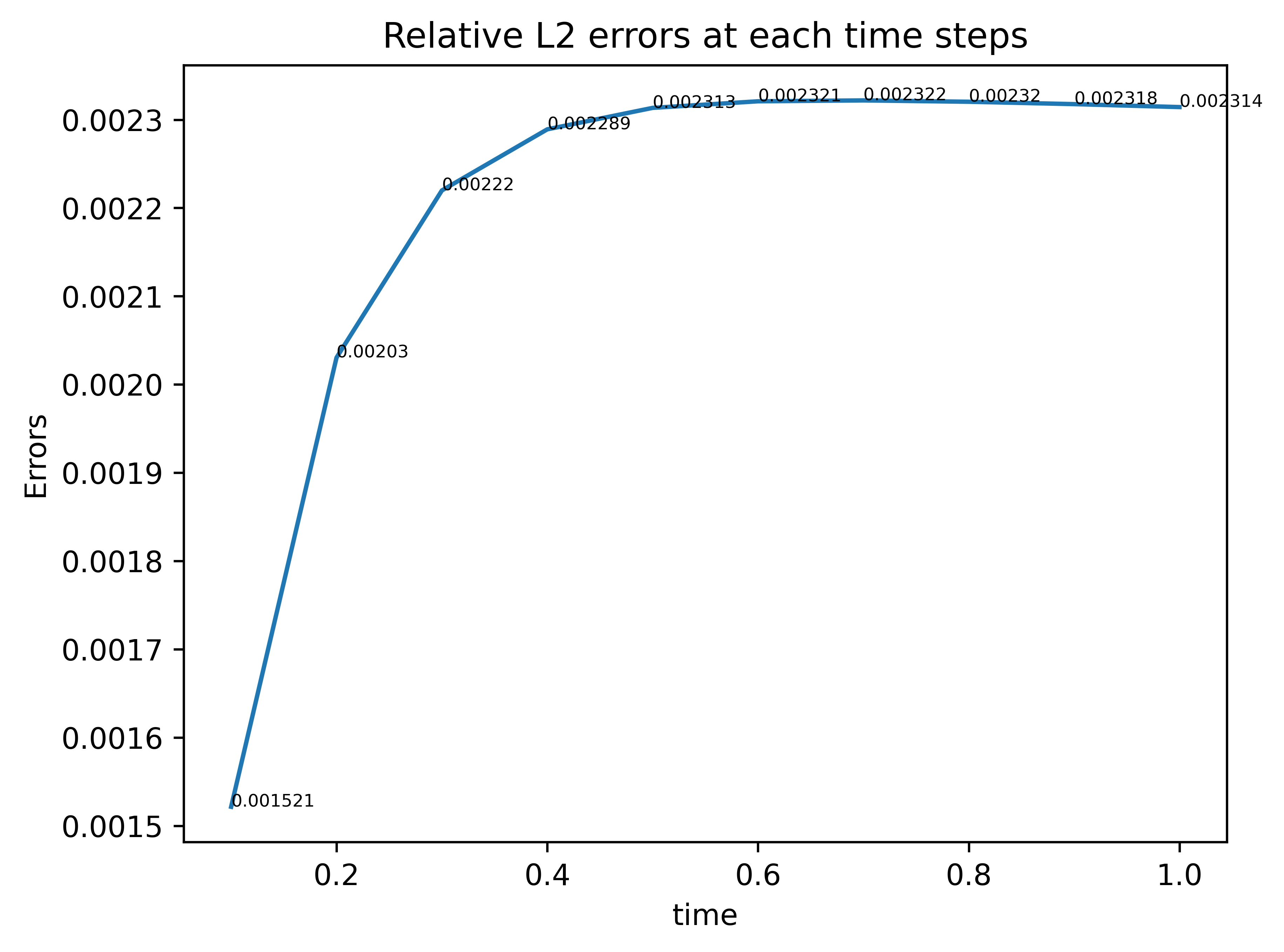}
        \label{fig1127_4_1_b}
    \end{subfigure}
    \caption{Numerical results for heat equation \eqref{eq6} by backward Euler scheme with $\tau = 0.1$. The second graph illustrates the error at each time step, indicating that it increases at a relatively slow rate.}
    \label{fig1127_4_1}
\end{figure}

\begin{table}[htbp]
\centering
\begin{tabular}{|c|c|p{2cm}|p{2cm}|p{2cm}|p{2cm}|}
\hline
$a$ & $b$ & Absolute $L_2$ error&  Absolute $L_{\infty}$ error&  Relative $L_2$  error& Relative $L_{\infty}$ error \\
\hline
0.5591 & 0.8291 & 0.000194 & 0.000819 & 0.0808\% & 0.2265\% \\
0.4406 & 0.8977 & 0.000191 & 0.000793 & 0.0817\% & 0.2216\% \\
0.3153 & 0.9490 & 0.000188 & 0.000755 & 0.0834\% & 0.2153\% \\
0.5571 & 0.8305 & 0.000193 & 0.000819 & 0.0808\% & 0.2264\% \\
\hline
\end{tabular}
\caption{Numerical errors for different values of \(a\) and \(b\).}
\label{lab0220_1}
\end{table}
\item \textbf{Solving multiple PDEs.} In addition to the backward Euler scheme, the Crank–Nicolson scheme described in Section~\ref{Heat Equations} is employed to solve the heat equation. Leveraging the pre-trained \textbf{Model 1} and \textbf{Model 3} from Section \ref{Modified Helmholtz Equations}, we are able to efficiently solve a large number of heat equations of the form \eqref{eq6} with coefficients \( a \sim \mathcal{U}(0.25, 0.75) \) and \( b = \sqrt{1 - a^2} \).Under the specified boundary and initial conditions, and after precomputing the partial derivatives of the Green's function (which can be stored for repeated use), a total of 10,000 equations with randomly sampled values of $a$ and $b$ are solved simultaneously in just about \textbf{0.0975} seconds using a single RTX 4090 GPU. Numerical errors at final time \( T = 1 \), using a time step size of \( \tau = 0.1 \), are reported in Table~\ref{lab0220_1}. As the table shows, the proposed model achieves consistently high accuracy across the sampled equations while demonstrating exceptional computational efficiency.

\item \textbf{Uncertain quantification.} To demonstrate the potential of our operator learning method in solving random PDEs, we consider the heat equation \eqref{eq6} with uncertainty in the parameter $a$. Specifically, assume that
\[
a \sim \mathcal{N}(0.5, 0.05^2),
\]
and apply the truncation $a = \texttt{np.clip}(a, 0.2, 0.8)$ to ensure $a$ remains in the admissible range. The coefficient $b$ is then computed as $b = \sqrt{1 - a^2}$ to ensure the identity $a^2 + b^2 = 1$ still holds. Using the Crank-Nicolson scheme described in Section \ref{Heat Equations}, together with the pre-trained \textbf{Model 1} and \textbf{Model 3} from Section \ref{Modified Helmholtz Equations}, we solve \eqref{eq6} for $M = 10,000$ different realizations of $a$. Figure~\ref{fig:UQ_histograms} shows the distributions of: (a) the sampled values of $a$, (b) the exact solution $u(x, t)$ at point $(0.43, 0.2)$ and $t = 1$, (c) the predicted value by the model at the same location and time, and (d) the pointwise absolute error. From the histograms, it is evident that the method introduced in this paper provides not only highly accurate predictions, but also captures the statistical distribution of the output under input uncertainty with great fidelity. The predicted solutions not only match the exact values pointwise, but also preserve key statistical properties. Additionally, in the global sense, the empirical mean and standard deviation of the predicted values closely match those of the exact solutions, as shown in Table \ref{tab0619}.
\begin{figure}[htbp]
    \centering
    \includegraphics[width=0.9\textwidth]{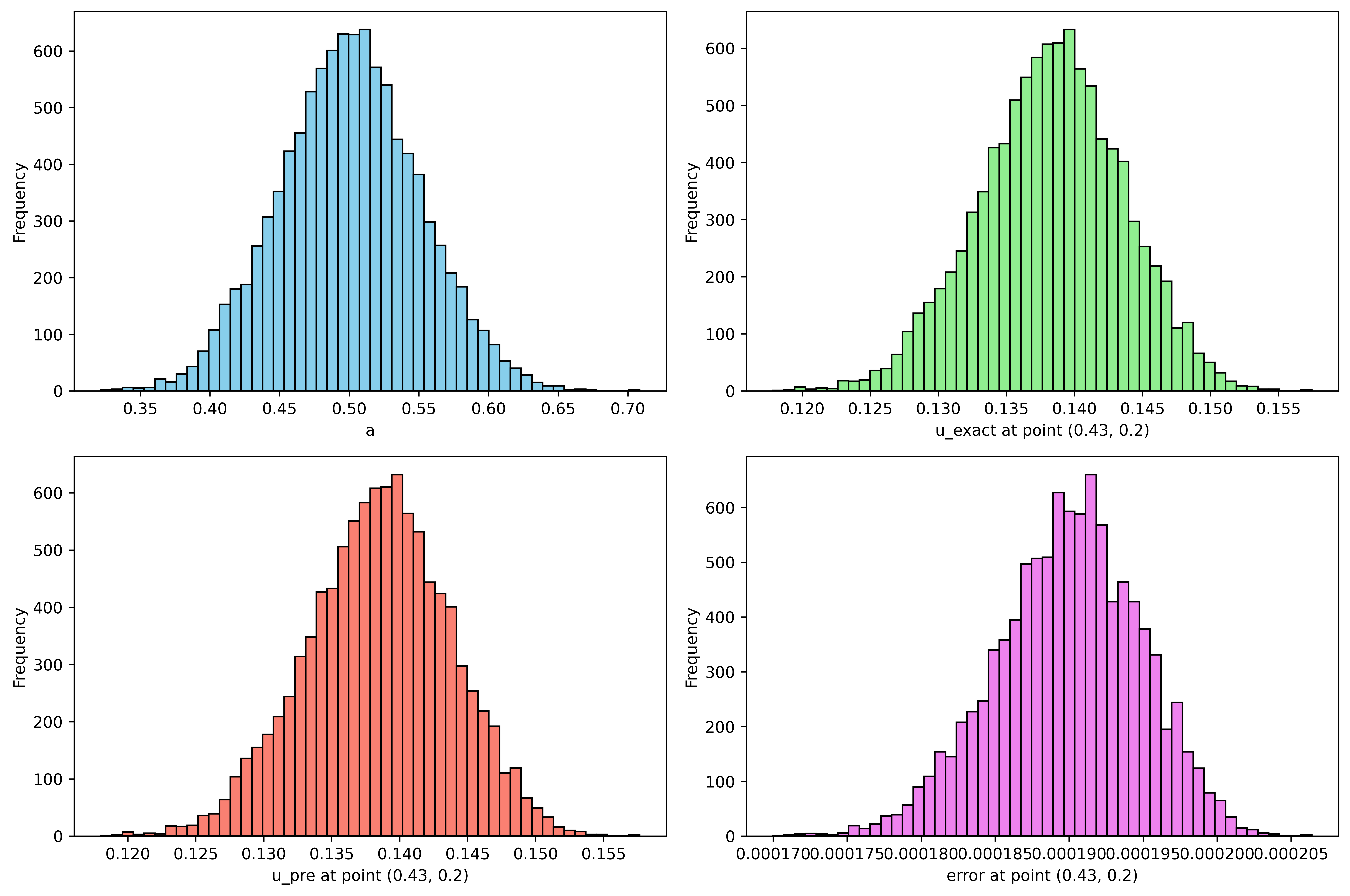}
    \caption{UQ results for the heat equation with uncertain parameter $a$: 
    (Top Left) Histogram of sampled $a$; 
    (Top Right) Exact solution $u$ at $(0.43, 0.2)$; 
    (Bottom Left) Predicted $u$ at $(0.43, 0.2)$; 
    (Bottom Right) Absolute error between predicted and exact values.}
    \label{fig:UQ_histograms}
\end{figure}

\begin{table}[htbp]
\centering
\begin{tabular}{l r}
\toprule
\textbf{Metric} & \textbf{Value} \\
\midrule
Mean of Exact Solution      & 0.221880 \\
Std of Exact Solution       & 0.082560 \\
Mean of Prediction          & 0.221832 \\
Std of Prediction           & 0.082477 \\
Mean of Error                  & -0.000047 \\
Std of Error                & 0.000186 \\
Max Absolute Error          & 0.000838 \\
Relative L2 Error  & \textbf{0.000812} \\
95\% Quantile Error         & 0.000404 \\
\bottomrule
\end{tabular}
\caption{Global statistics computed over all $10{,}000$ samples and $1{,}681$ grid point.}
\label{tab0619}
\end{table}

The global error analysis across all 10,000 samples and 1,681 grid points demonstrates exceptional predictive accuracy. The mean prediction ($0.221832$) is nearly identical to the mean exact solution ($0.221880$). The standard deviation of the prediction ($0.082477$) also closely matches that of the exact solution ($0.082560$), indicating excellent preservation of the solution's variability. Most notably, the relative L2 error of $\mathbf{8.12 \times 10^{-4}}$ — less than one-tenth of a percent — signifies that the model's predictions are extremely close to the ground truth in an integrated (L2) sense. Furthermore, the 95th percentile error of $4.04 \times 10^{-4}$ confirms that 95\% of all predicted values across all samples and grid points deviate from the true values by less than 0.0004, which is remarkably small for most scientific or engineering applications. This level of precision strongly suggests that the model has successfully learned the underlying function space and generalizes well across the entire domain.

\end{itemize}

\subsection{Wave Equations} \label{Numerical Wave Equations}
Consider the classical wave equation on a two-dimensional square domain $\Omega = [0, 1]^2$ with time-dependent Dirichlet boundary conditions:
\begin{equation} \label{eq0220_1}
\begin{cases} 
u_{tt} = \Delta u , & \mathbf{x} \in \Omega,\ t \in [0, T], \\ 
u(\mathbf{x},0) = \sin(a x_1 + b x_2) , & \mathbf{x} = (x_1,x_2) \in \Omega, \\ 
u_t(\mathbf{x},0) = -\cos(a x_1 + b x_2) , & \mathbf{x} \in \Omega, \\
u(\mathbf{x},t) = \sin(a x_1 + b x_2 - t), & \mathbf{x} \in \partial \Omega,\ t \in [0, T],
\end{cases}
\end{equation}
where $a, b \in (0, 1)$ satisfy $a^2 + b^2 = 1$. The exact solution is given by
\[
u(\mathbf{x}, t) = \sin(a x_1 + b x_2 - t).
\]

\begin{figure}[htbp]
    \centering
    \includegraphics[width=0.82\textwidth]{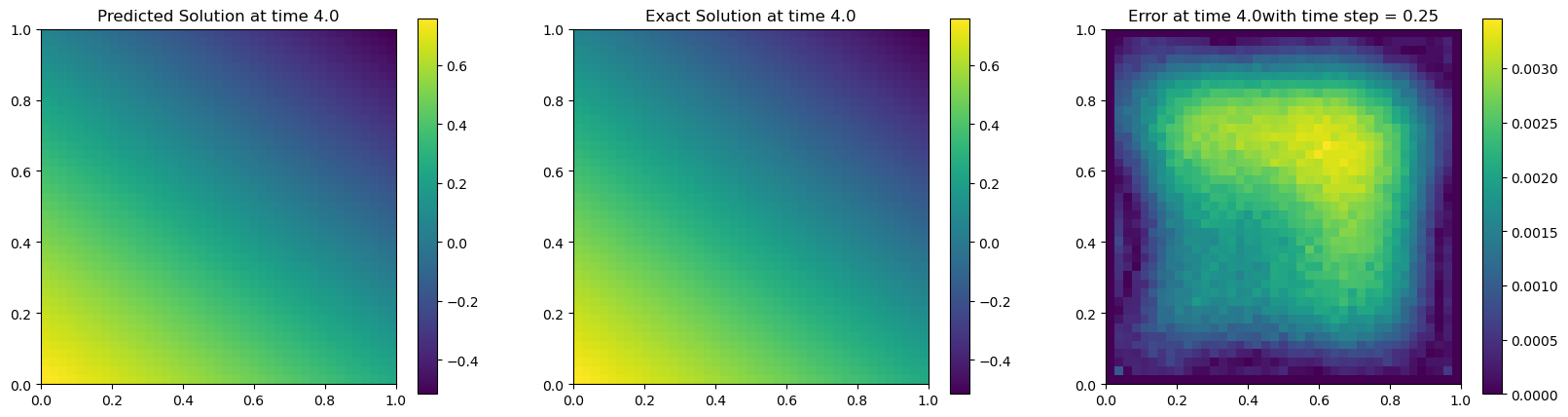} 
    \caption{Numerical results at time $T = 4$ for wave equation \eqref{eq0220_2} with $a = 0.6, \tau = \frac{1}{4}$.}
    \label{fig0220_4}
\end{figure}

\begin{figure}[htbp]
    \centering
    \includegraphics[width=0.82\textwidth]{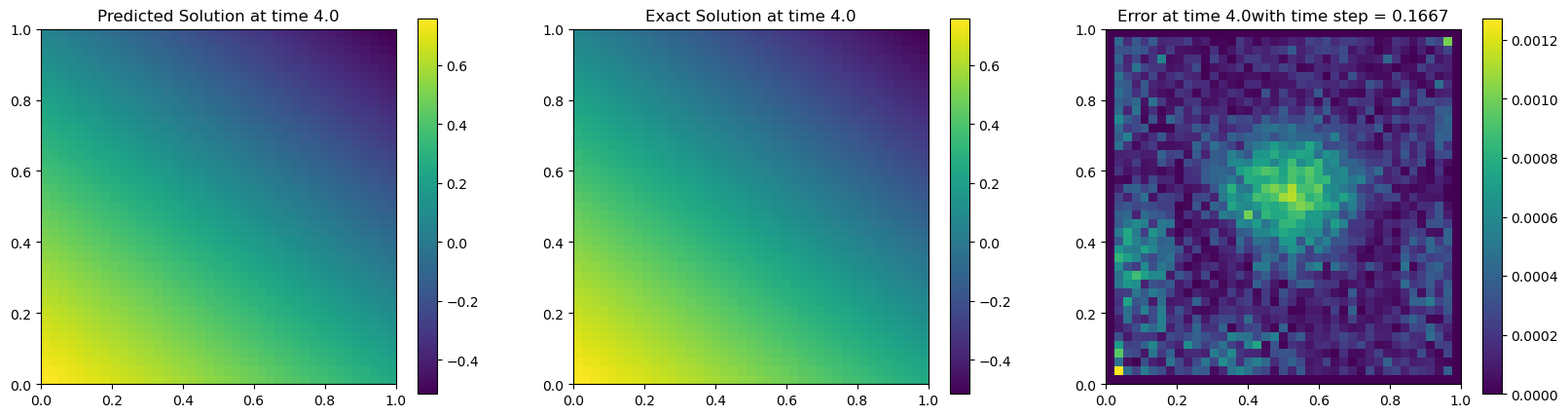} 
    \caption{Numerical results at time $T = 4$ for wave equation \eqref{eq0220_2} with $a = 0.6, \tau = \frac{1}{6}$.}
    \label{fig0220_3}
\end{figure}

\begin{figure}[htbp]
    \centering
    \includegraphics[width=0.82\textwidth]{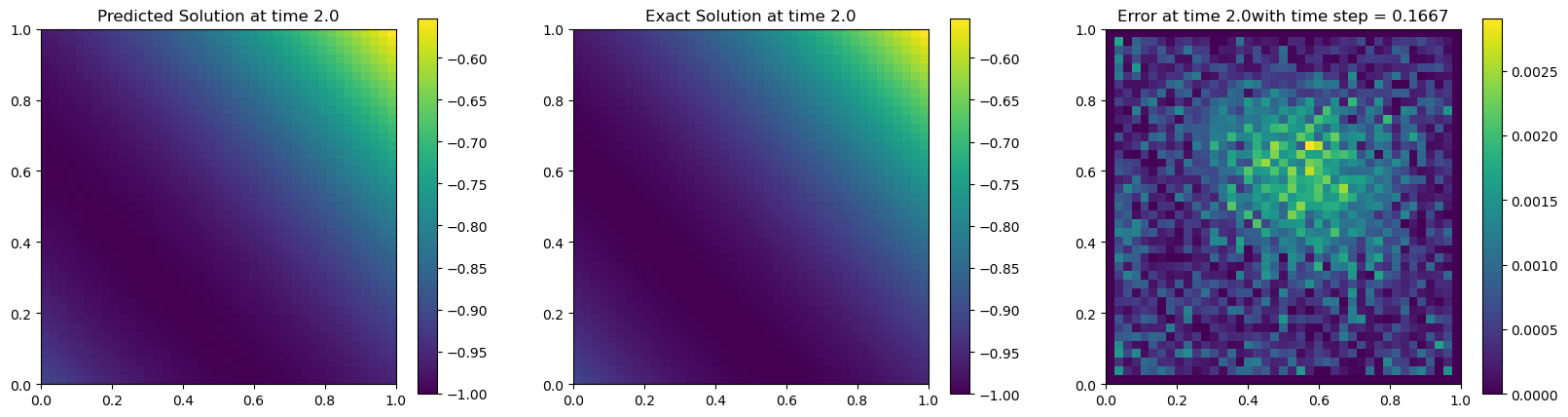} 
    \caption{Numerical results at time $T = 2$ for wave equation \eqref{eq0220_2} with $a = 0.7335, \tau = \frac{1}{6}$.}
    \label{fig0220_5}
\end{figure}

To discretize the equation in time, we apply the $\theta$-scheme with $\theta = \frac{1}{2}$ introduced in Section \ref{Wave Equations}, and the resulting scheme reads:
\begin{equation} \label{eq0220_2}
\begin{cases}
(\mathcal{I} - \frac{\tau^2}{2} \Delta) u^{n+1} = 2 u^n - u^{n-1} + \frac{\tau^2}{2} \Delta u^{n-1} := F^{n+1}, & \mathbf{x} \in \Omega, \\
u^{n+1} = \sin(a x_1 + b x_2 - (n + 1)\tau), & \mathbf{x} \in \partial \Omega.
\end{cases}
\end{equation}
Here, $\mathcal{I}$ denotes the identity operator. As in the heat equation case, equation \eqref{eq0220_2} at each time step is an elliptic PDE with Dirichlet boundary conditions. Therefore, the above wave equation can be efficiently solved using the pre-trained \textbf{Model 2} and \textbf{Model 4} described in Section~\ref{Modified Helmholtz Equations}. These models approximate the solution of elliptic PDEs of the form $(\mathcal{I} - \lambda \Delta)u = f$ with high accuracy. Representative numerical results illustrating the performance of our method are presented in Table \ref{tab0220_1}, Figure \ref{fig0220_4}, Figure \ref{fig0220_3}, and Figure \ref{fig0220_5}.

\begin{table}[htbp]
    \centering
    \begin{tabular}{|c|c|c|p{2cm}|p{2cm}|p{2cm}|p{2cm}|} \hline 
         $\tau$ & a & $T$ & Absolute $L_2$ error&  Absolute $L_{\infty}$ error&  Relative $L_2$  error& Relative $L_{\infty}$ error\\ \hline 
         $\frac{1}{4}$ & 0.6 & 4 & 0.001751 & 0.003456 & 0.5469\% & 0.4566\% \\ \hline 
         $\frac{1}{6}$ & 0.6 & 4 & 0.000320 & 0.001271 & 0.1001\% & 0.1679\% \\ \hline 
         $\frac{1}{6}$ & 0.7335 & 2 & 0.000831 & 0.002900 & 0.0899\% & 0.2900\% \\ \hline 
    \end{tabular}
    \caption{Numerical errors for wave equation \eqref{eq0220_2} solved by $\theta-$scheme and NEKM under different coefficient $a$, time step $\tau$ and end time $T$.}
    \label{tab0220_1}
\end{table}

\subsection{Schrödinger equations} \label{Schrödinger equations}
As introduced in Section \ref{The Schrödinger equation}, whether using the Strang splitting or Lie-Trotter splitting methods to solve the Schrödinger equations, solving the linear part of the PDE—namely equation \eqref{eq0608_S3}—is crucial. After applying the backward Euler time discretization, the problem reduces to solving PDEs of the type given in the third equation of \eqref{eq0608_S6} or the second equation of \eqref{eq0608_S8}. We emphasize that, since \( u \) is a complex-valued function, the resulting equations (such as the third equation in \eqref{eq0608_S6} or the second equation in \eqref{eq0608_S8}) can be reformulated as a system of PDEs for the real and imaginary parts of \( u \). As stated in Appendix \ref{Operator Learning for PDE Systems}, the PDE system of interest at this stage is
\begin{equation} \label{eq1124_11}
    \begin{cases} \mathbf{L}_\lambda \mathbf{u} = \mathbf{f}, & \mathbf{x} \in \Omega \\ 
    \mathbf{u} = \mathbf{g}^D, & \mathbf{x} \in \partial \Omega, \end{cases}
\end{equation}
in which $\mathbf{u} = (u_1, u_2)^\top, \mathbf{f} = (f_1, f_2)^\top, \mathbf{g}^D = (g_1, g_2)^\top$ and 
\begin{equation} \label{eq1124_111}
\mathbf{L}_\lambda = \begin{pmatrix}
I & -\lambda \Delta \\
\lambda \Delta & I
\end{pmatrix}
\end{equation}
is a matrix of operators. The boundary integral formulation of such systems, along with the implementation of operator learning (which, although similar to the process described in Section \ref{Operator Learning}, includes several necessary adjustments), is provided in Appendix \ref{Operator Learning for PDE Systems}. Here, we only present the corresponding numerical results.

\begin{itemize} 
\item \textbf{Operator learning for source terms in systems.} Consider the PDE system \eqref{eq1124_11} on domain $\Omega = [0, 1]^2$ with homogeneous Dirichlet boundary conditions, i.e., $g^D_\text{Re} = g^D_\text{Im} \equiv 0$. Following the operator learning framework described in Section \ref{Operator Learning for Source Terms}, we train \textbf{Model 5} $\mathcal{NN}_1(\texttt{f}, \texttt{k}, \texttt{x})$ for $\lambda \in [0.05, 0.1]$. For each model, 11 uniformly sampled values of $\lambda$ are used for training. The source terms $f_\text{Re}$ and $f_\text{Im}$ are generated using a similar strategy as in Section~\ref{Modified Helmholtz Equations}. The test solutions are chosen as follows:
\[
u_\text{Re}(\mathbf{x}) = \sin(\pi x_1)\, x_2 (1 - x_2), \quad
u_\text{Im}(\mathbf{x}) = x_1 (1 - x_1)\, \sin(\pi x_2),
\]
where $\mathbf{x} = (x_1, x_2) \in \Omega$. Some representative results for \textbf{Model 5} are shown in Table \ref{tab0311_4}.
\begin{table}[htbp]
    \centering
    \renewcommand{\arraystretch}{1.2}
    \begin{tabular}{c|c|c|c}
        \hline
        $\lambda$ & Error Type & $u_1$ & $u_2$ \\
        \hline
        \multirow{4}{*}{0.1} 
        & Absolute $L_2$ Norm Error & 0.00027043 & 0.00036826 \\
        & Absolute Max Norm Error & 0.00171460 & 0.00179922 \\
        & Relative $L_2$ Norm Error (\%) & 0.21471 & 0.29239 \\
        & Relative Max Norm Error (\%) & 0.68584 & 0.71969 \\
        \hline
        \multirow{4}{*}{$\frac{1}{15}$} 
        & Absolute $L_2$ Norm Error & 0.00053998 & 0.00036784 \\
        & Absolute Max Norm Error & 0.00151578 & 0.00109008 \\
        & Relative $L_2$ Norm Error (\%) & 0.42872 & 0.29206 \\
        & Relative Max Norm Error (\%) & 0.60631 & 0.43603 \\
        \hline
        \multirow{4}{*}{0.05} 
        & Absolute $L_2$ Norm Error & 0.00037180 & 0.00035480 \\
        & Absolute Max Norm Error & 0.00161524 & 0.00149200 \\
        & Relative $L_2$ Norm Error (\%) & 0.29520 & 0.28170 \\
        & Relative Max Norm Error (\%) & 0.64609 & 0.59680 \\
        \hline
    \end{tabular}
    \caption{Numerical results for \textbf{Model 5}. Note that the value $\lambda = \frac{1}{15}$ is not included in the training set $\{\lambda_i \}_{i = 1}^{11}$, and the test functions $u_\text{Re}$ and $u_\text{Im}$ are not included in the training data.}
    \label{tab0311_4}
\end{table}

\item \textbf{Operator learning for boundary conditions in systems.} Consider the PDE system \eqref{eq1124_11} on the domain $\Omega = [0, 1]^2$ with vanishing source terms, i.e., $f_\text{Re} = f_\text{Im} \equiv 0$. Based on the boundary integral network framework introduced in Section \ref{Operator Learning for Boundary Conditions}, we train \textbf{Model 6} $\mathcal{NN}_2(\texttt{k}, \texttt{g}^\texttt{D})$ for $\lambda \in [0.05, 0.1]$. For each model, 11 uniformly sampled values of $\lambda$ are used for training. The boundary data $g^D_\text{Re}$ and $g^D_\text{Im}$ are generated using a similar strategy as in Section \ref{Modified Helmholtz Equations}. The test solutions used in this section are given by
\[
u_\text{Re}(\mathbf{x}) = G_{\lambda,0}^{1,1}(\mathbf{x}, \mathbf{y}), \quad u_\text{Im}(\mathbf{x}) = G_{\lambda,0}^{2,1}(\mathbf{x}, \mathbf{y}),
\]
where $\mathbf{y} = (1.2, 1.2) \in \mathbb{R}^2$. The functions $G_{\lambda,0}^{1,1}$ and $G_{\lambda,0}^{2,1}$ are the fundamental solutions defined in equation \eqref{eq0310_5}, with their explicit expressions provided in Appendix \ref{Boundary Integral Derivation for PDE System}. Representative results for \textbf{Model 6} are reported in Table \ref{tab0311_2}.
\begin{table}[htbp]
    \centering
    \renewcommand{\arraystretch}{1.2}
    \begin{tabular}{c|c|c|c}
        \hline
        $\lambda$ & Error Type & $u_1$ & $u_2$ \\
        \hline
        \multirow{4}{*}{0.1} 
        & $L_\infty$ Absolute Error & 0.00048173 & 0.00034244 \\
        & $L_2$ Absolute Error & 0.00023093 & 0.00012468 \\
        & $L_\infty$ Relative Error & 0.00066680 & 0.00101281 \\
        & $L_2$ Relative Error & 0.00129046 & 0.00143291 \\
        \hline
        \multirow{4}{*}{$\frac{1}{15}$} 
        & $L_\infty$ Absolute Error & 0.00056499 & 0.00029418 \\
        & $L_2$ Absolute Error & 0.00010844 & 9.8777e-05 \\
        & $L_\infty$ Relative Error & 0.00063751 & 0.00128101 \\
        & $L_2$ Relative Error & 0.00062284 & 0.00100027 \\
        \hline
        \multirow{4}{*}{0.05} 
        & $L_\infty$ Absolute Error & 0.00054753 & 0.00037954 \\
        & $L_2$ Absolute Error & 0.00012341 & 0.00015064 \\
        & $L_\infty$ Relative Error & 0.00055868 & 0.00168012 \\
        & $L_2$ Relative Error & 0.00075427 & 0.00140017 \\
        \hline
    \end{tabular}
    \caption{Numerical results for \textbf{Model 6}. Note that the value $\lambda = \frac{1}{15}$ is not included in the training set $\{\lambda_i \}_{i = 1}^{11}$, and the test functions $u_\text{Re}$ and $u_\text{Im}$ are not included in the training data.}
    \label{tab0311_2}
\end{table}

\item \textbf{Solving the Schrödinger equation.}
In the following, consider a specific Schrödinger equation given by equation~\eqref{eq0608_S1}, where the potential function is defined as \( v(\mathbf{x}) = v(x_1, x_2) \\ = 1 - \cos^2(x_1) \cos^2(x_2) \), with \( w = 1 \) and the spatial domain \( \Omega = [0, 1]^2 \). The initial and boundary conditions are specified based on the exact solution \( u(\mathbf{x}, t) = e^{it} \cos(x_1) \cos(x_2) \). To solve this equation, we employ the Strang splitting method~\eqref{eq0608_S6} and the Lie–Trotter splitting method~\eqref{eq0608_S8}. By separating the real and imaginary parts of the PDEs, the coupled system~\eqref{eq0310_4} is obtained, which is studied in Appendix~\ref{Operator Learning for PDE Systems}. This allows to directly apply the previously trained operator learning \textbf{Models 5} and \textbf{6}. Numerical results for Lie splitting and Strang splitting are presented in Table \ref{tab1124_1} and Table \ref{tab1124_2}, respectively. It is worth noting that, in this example, the error is computed based on the solutions at all time steps. The evolution diagram of the numerical results calculated by Strang splitting with $\tau = 0.12$ are shown in Figure \ref{fig1124_1}. These results demonstrating that NEKM can be effectively integrated with various operator splitting schemes while maintaining both stability and high accuracy in solving time-dependent Schrödinger equations.

\begin{table}[htbp]
\centering
\begin{tabular}{c|c|ccc}
\hline
Quantity & Error Type & $\tau = 0.08$ & $\tau = 0.07$ & $\tau = 0.06$ \\
\hline
\multirow{4}{*}{$u_1$ (Real part)}
& Absolute L2   & 0.004171 & 0.003174 & 0.002559 \\
& Absolute Max  & 0.011760 & 0.009386 & 0.007782 \\
& Relative L2   & 0.008426 & 0.005980 & 0.004485 \\
& Relative Max  & 0.011798 & 0.009409 & 0.007796 \\
\hline
\multirow{4}{*}{$u_2$ (Imaginary part)}
& Absolute L2   & 0.004086 & 0.003284 & 0.003017 \\
& Absolute Max  & 0.012876 & 0.009785 & 0.008988 \\
& Relative L2   & 0.007708 & 0.006643 & 0.006737 \\
& Relative Max  & 0.012882 & 0.009929 & 0.009643 \\
\hline
\multirow{4}{*}{$(u_1,u_2)$ Combined}
& Absolute L2   & 0.004129 & 0.003229 & 0.002797 \\
& Absolute Max  & 0.012876 & 0.009785 & 0.008988 \\
& Relative L2   & 0.008051 & 0.006297 & 0.005454 \\
& Relative Max  & 0.012882 & 0.009809 & 0.009004 \\
\hline
\end{tabular}
\caption{Global errors over all time steps for the Schr\"odinger equation using Lie splitting with different time step sizes $\tau$. The total number of time steps is $N = 20$.}
\label{tab1124_1}
\end{table}

\begin{table}[htbp]
\centering
\begin{tabular}{c|c|ccc}
\hline
Quantity & Error Type & $\tau = 0.16$ & $\tau = 0.14$ & $\tau = 0.12$ \\
\hline
\multirow{4}{*}{$u_1$ (Real part)}
& Absolute L2   & 0.004449 & 0.003039 & 0.002492 \\
& Absolute Max  & 0.014634 & 0.010459 & 0.010305 \\
& Relative L2   & 0.009236 & 0.005863 & 0.004449 \\
& Relative Max  & 0.014823 & 0.010562 & 0.010380 \\
\hline
\multirow{4}{*}{$u_2$ (Imaginary part)}
& Absolute L2   & 0.004032 & 0.002530 & 0.002121 \\
& Absolute Max  & 0.013146 & 0.008458 & 0.007795 \\
& Relative L2   & 0.007435 & 0.004987 & 0.004603 \\
& Relative Max  & 0.013152 & 0.008583 & 0.008363 \\
\hline
\multirow{4}{*}{$(u_1,u_2)$ Combined}
& Absolute L2   & 0.004245 & 0.002796 & 0.002314 \\
& Absolute Max  & 0.014634 & 0.010459 & 0.010305 \\
& Relative L2   & 0.008278 & 0.005452 & 0.004512 \\
& Relative Max  & 0.014640 & 0.010562 & 0.010380 \\
\hline
\end{tabular}
\caption{Global errors over all time steps for the Schr\"odinger equation using Strang splitting with different time step sizes $\tau$. The total number of time steps is $N = 10$.}
\label{tab1124_2}
\end{table}

\begin{figure}[htbp]
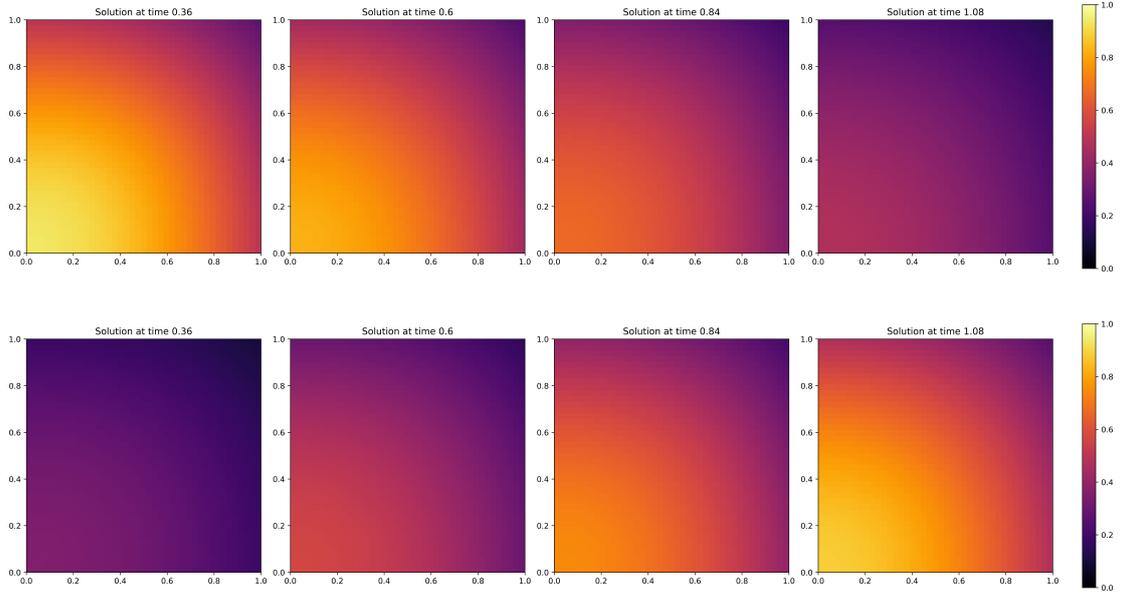

    \centering
    \begin{subfigure}
        \centering
        \includegraphics[width=\textwidth]{figures/fig1124_1.png}
        \label{fig:subfig1}
    \end{subfigure}
    \begin{subfigure}
        \centering
        \includegraphics[width=\textwidth]{figures/fig1124_2.png}
        \label{fig:subfig2}
    \end{subfigure}
    \caption{The numerical solutions to the Schrödinger equation using Strang splitting. The first row displays the real parts, while the second row shows the imaginary parts.}
    \label{fig1124_1}
\end{figure}

\end{itemize}

\subsection{Heat Equations in a Petal-Shaped Domain} \label{Heat Equations in a Petal-Shaped Domain}
\begin{figure}[htbp]
    \centering
    \includegraphics[width=0.6\textwidth]{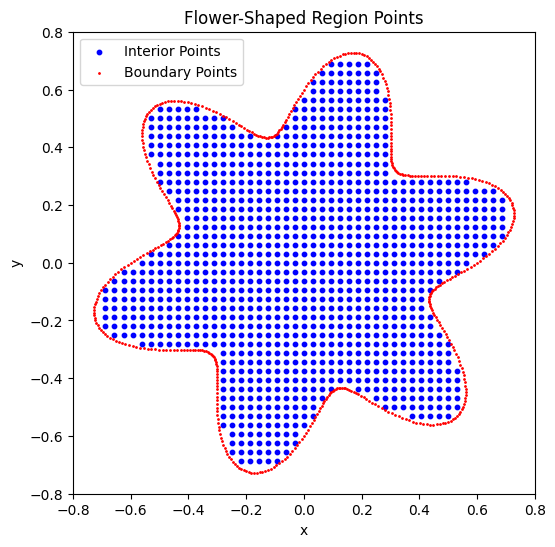} 
    \caption{Schematic diagram of the petal-shaped domain.}
    \label{fig0422_1}
\end{figure}

\begin{figure}[htbp]
    \centering
    \begin{subfigure}
        \centering
        \includegraphics[width=\textwidth]{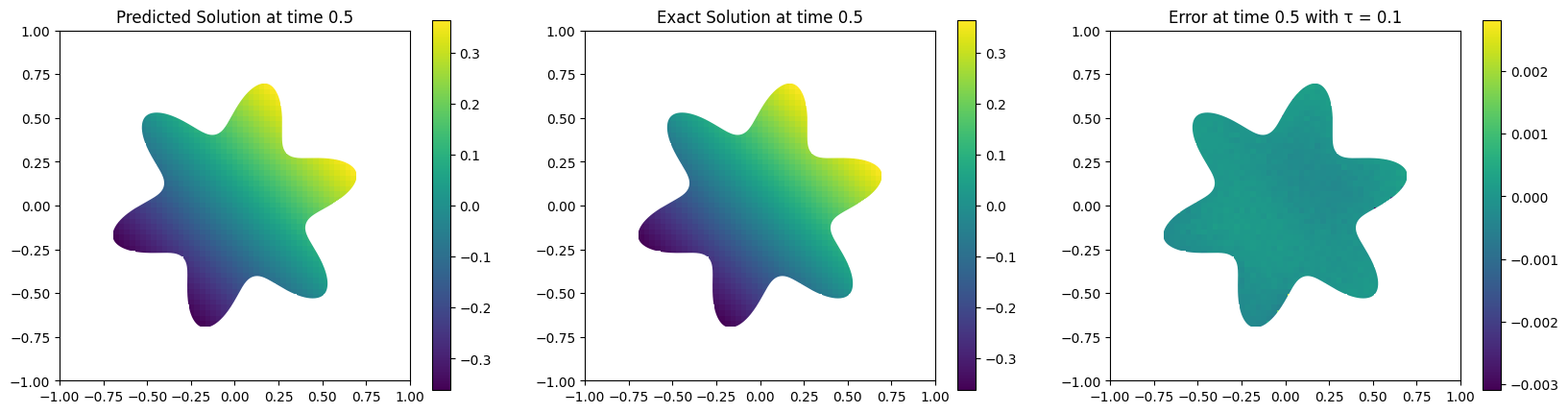}
    \end{subfigure}
    \vspace{0.5cm} 
    \begin{subfigure}
        \centering
        \includegraphics[width=\textwidth]{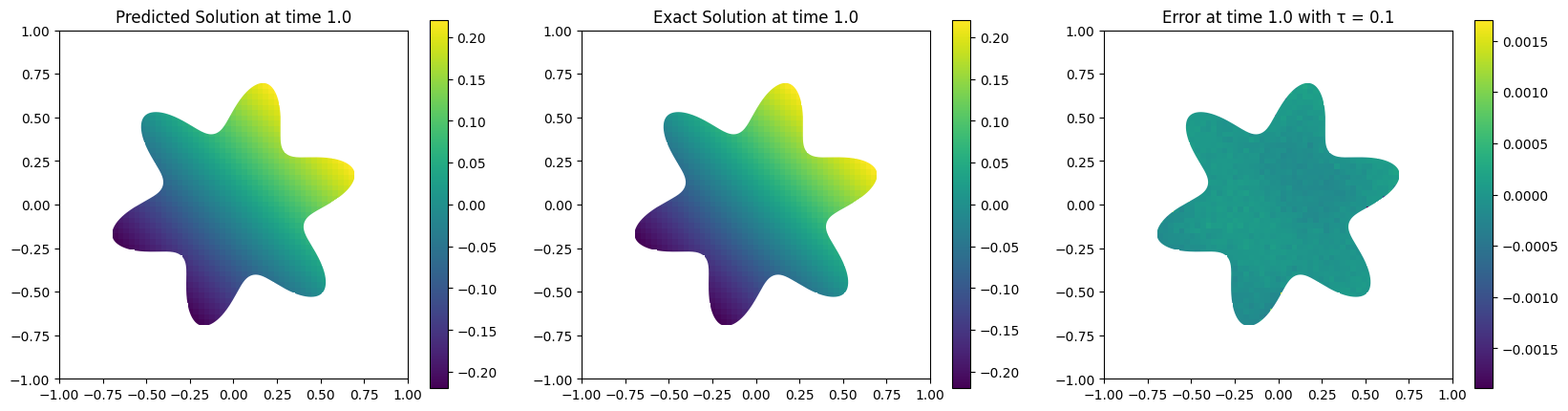}
    \end{subfigure}
    \vspace{0.5cm} 
    \begin{subfigure}
        \centering
        \includegraphics[width=\textwidth]{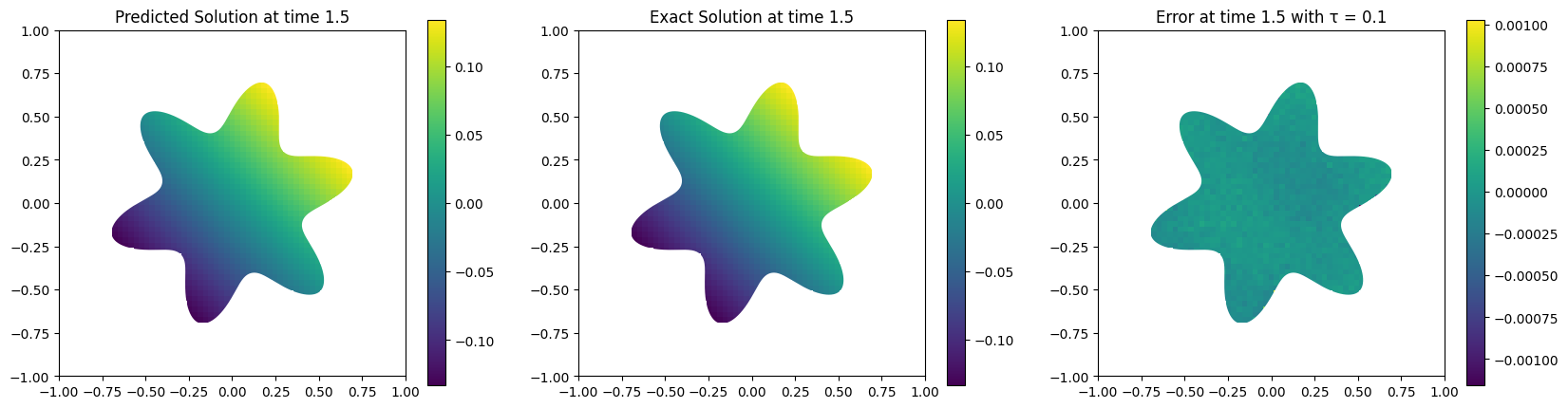}
    \end{subfigure}
    \vspace{0.5cm} 
    \begin{subfigure}
        \centering
        \includegraphics[width=\textwidth]{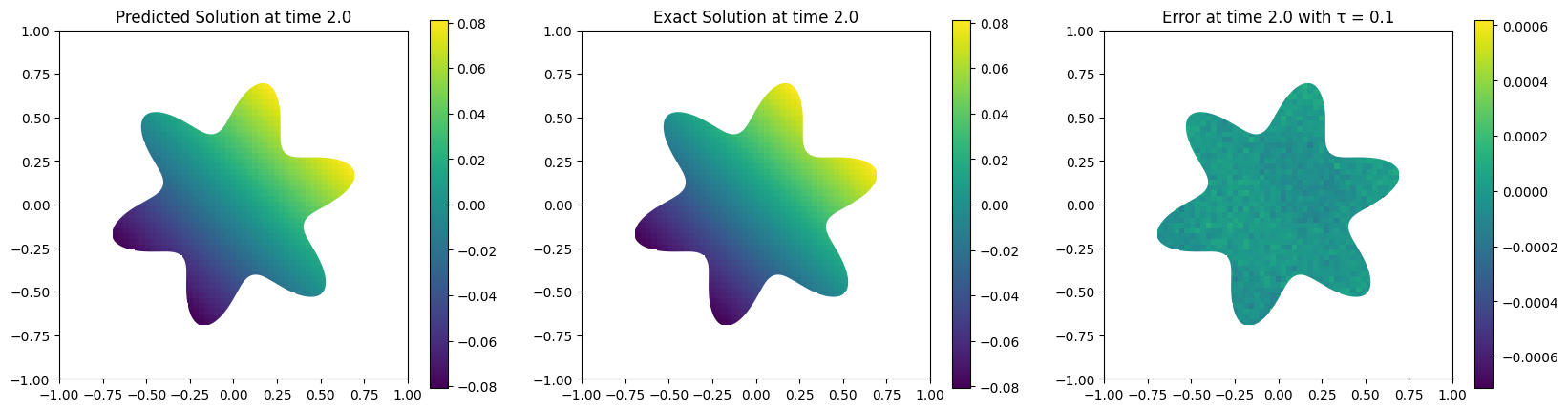}
    \end{subfigure}
    \caption{Some numerical solutions for the heat equation in the petal-shaped domain with $\tau = 0.1$.}
    \label{fig0422_2}
\end{figure}

In this section, we change the computational domain from the previously considered rectangular region to a more complex petal-shaped domain (illustrated in Figure~\ref{fig0422_1}) to further validate the effectiveness and generalization ability of NEKM on irregular geometries. The parametric equation of the petal-shaped boundary is given by 
\[
\{(x(\theta), y(\theta))\ |\ x(\theta) = 0.6(1 + 0.25 \sin(6\theta)) \cos(\theta),\ y(\theta) = 0.6(1 + 0.25 \sin(6\theta)) \sin(\theta)\}.
\]
Consider the heat equation \( u_t = \Delta u \) on this domain, with both the initial and Dirichlet boundary conditions derived from the exact solution
\[
u(\mathbf{x}, t) = \exp(-t) \sin\left(\tfrac{\sqrt{2}}{2} x_1\right) \sin\left(\tfrac{\sqrt{2}}{2} x_2\right).
\]
The numerical error is evaluated using over 1,000 uniformly distributed interior test points within the domain. Other experimental settings remain consistent with those in Sections~\ref{Modified Helmholtz Equations} and~\ref{Numerical Heat Equations}. Representative results are presented in Figure~\ref{fig0422_2} and Table~\ref{tab0422}. These results confirm that NEKM framework maintains robustness and accuracy even on complex, non-rectangular domains. The method exhibits strong generalization ability, enabling stable and high-precision solutions for time-evolutionary partial differential equations beyond standard geometries.

\begin{table}[htbp]
\centering
\begin{tabular}{c|c|c|c|c}
\hline
Time $T$ & Absolute $L^2$ Error & Absolute Max Error & Relative $L^2$ Error & Relative Max Error \\
\hline
0.5 & 0.000131 & 0.003097 & 0.1357\% & 0.8540\% \\
1.0 & 0.000081 & 0.001891 & 0.1382\% & 0.8596\% \\
1.5 & 0.000051 & 0.001159 & 0.1445\% & 0.8690\% \\
2.0 & 0.000035 & 0.000716 & 0.1606\% & 0.8843\% \\
\hline
\end{tabular}
\caption{Error metrics for the heat equation on the petal-shaped domain with $\tau = 0.1$ at different time levels}
\label{tab0422}
\end{table}

\section{Conclusion} \label{Conclusion}
This work proposes the \textit{Neural Evolutionary Kernel Method (NEKM)}, a novel operator learning framework for solving a broad class of time-dependent partial differential equations (PDEs). By integrating the boundary integral method (BIM) into the architecture of deep neural networks, NEKM leverages the mathematical structure of elliptic operators to construct efficient and accurate solvers for time-dependent problems. Specifically, the solution of elliptic PDEs is decomposed into source-driven and boundary-driven components, and two neural networks are trained to approximate the corresponding solution operators. Through suitable time discretization schemes—such as backward Euler, Crank--Nicolson, and operator splitting—the framework extends naturally to time-dependent PDEs, including the heat equation, wave equation, and Schr\"{o}dinger equation.

A key strength of NEKM lies in its \textit{generality and scalability}. Since the models operate at the operator level and are parameterized by PDE coefficients, they generalize effectively across families of equations, boundary conditions, and parameter regimes. Numerical experiments on both classical rectangular domains and more challenging petal-shaped geometries confirm that NEKM achieves high accuracy and stability. Moreover, the tensor-based implementation enables parallel prediction of thousands of PDE solutions, making the method particularly suitable for \textit{uncertainty quantification (UQ)} and ensemble simulations. Additionally, the boundary-driven component can be trained in a self-supervised manner, further improving data efficiency.

Looking ahead, the NEKM framework opens up several promising directions for future work:

\begin{itemize}
    \item \textbf{Application to phase field models.} Many phase field equations, such as the \textit{Allen--Cahn} and \textit{Cahn--Hilliard} equations, can be reformulated into sequences of elliptic subproblems using techniques such as Strang splitting or convex splitting \citep{li2022stability, guan2014second, wei2020integral}. This makes NEKM particularly suitable for efficiently solving such models, with potential applications in materials science, image processing, and interface dynamics.
    
    \item \textbf{Integration with inverse problems and uncertainty quantification.} As NEKM learns solution operators rather than individual solutions, it naturally supports fast evaluation under varying inputs. This makes it a powerful surrogate model for solving \textit{inverse problems}, conducting \textit{Bayesian inference}, or propagating input uncertainties in \textit{stochastic PDEs}. Future directions include coupling NEKM with variational inference, Hamiltonian Monte Carlo, or active learning to enable more robust and flexible modeling.
    
    \item \textbf{Adaptive time stepping for time-dependent PDEs.} Since the models allows for different PDE parameters as input, the method should be able to evolve time-dependent PDEs with different time steps at different moments, thereby implementing adaptive time stepping. However, for each time step, computing derivatives of the fundamental solution at all observation–integration point pairs can be expensive, especially for fundamental solutions involving special functions. A practical compromise is to allow the time step to be chosen from a small set of candidate values rather than continuously over a range, realizing a weak form of adaptive stepping without excessive computational cost.
    
    \item \textbf{Learning fundamental solutions for more general elliptic PDEs.} Currently, learning the mapping from boundary conditions to solutions of second-kind boundary integral equations requires explicit expressions of the fundamental solution. This restricts the application to more general elliptic PDEs, such as variable-coefficient problems. In future work, we aim to learn the fundamental solution of an elliptic operator itself. Concretely, for variable-coefficient elliptic PDEs, it may suffice to learn the matrix (or its inverse) formed by the integral kernel $K$ at all observation and integration points in \eqref{eq1119}. Such matrices are expected to have certain low-rank structures, which provides valuable insight for extending NEKM to a wider class of PDEs.
\end{itemize}

In summary, NEKM offers a knowledge-guided, efficient, and generalizable framework for solving time-evolution PDEs by combining operator learning with classical numerical structures, thereby providing a practical tool for scientific computing in solving time-dependent PDEs.

\bibliographystyle{plain} 
\bibliography{ref}  

\newpage
\section* {Appendix}
\appendix

\section{Fundamental Solution and Green's Function} \label{appendix A}
In this part, we make some supplementary introductions to the fundamental solution and Green's function. Let $\Omega \subset \mathbb{R}^d$ be a bounded domain, $\Omega^C:= \mathbb{R}^d \backslash \bar{\Omega}$ and $\Omega^* := \Omega$ or $\Omega ^ C$. $\mathcal{L}$ is an epplitic differential operator.

\begin{definition} \label{def A.1}
A function $\delta (x)$ is called a $d$-dimensional $\delta$-function if $\delta(x) \simeq 
\begin{cases}
0, & x \neq \vec 0, \\
\infty, &  x = \vec 0,
\end{cases}$ and which is also constrained to satisfy the identity $\int_{\mathbb{R}^d} \delta(x) dx = 0$. For all function $f$ that is continuous at $a \in \mathbb{R}^d$ we have $\int_{\mathbb{R}^d} f(x) \delta(x-a) dx = f(a)$.
\end{definition}

\begin{definition} \label{def A.2}
A function $G_0(x, y)$ is called the fundamental solution corresponding to equations $\mathcal{L} u (x)= 0$ if $G_0(x, y)$ is symmetric about $x$ and $y$ and $G_0(x, y)$ satisfy $\mathcal{L}_y G_0(x, y) = \delta(x - y)$, where $(x, y) \in \mathbb{R}^d \times \mathbb{R}^d$ and $\mathcal{L}_y$ is the differential operator $\mathcal{L}$ which acts on component y.
\end{definition}

\begin{definition} \label{def A.3}
A function $G(x, y)$ is called the Green's function corresponding to problem $$
\begin{cases}
\mathcal{L} u (x) = f(x)\ in\ \Omega^*, \\
u(x) = g(x)\ on\ \partial \Omega^*,
\end{cases}$$ if $G(x, y)$ is a $2d$-dimensional function satisfying $$
\begin{cases}
\mathcal{L}_y G(x, y) = \delta(x - y), & \forall x, y \in \Omega^*, \\
G(x, y) = 0, & \forall x \in \Omega^*, y \in \partial \Omega^*.
\end{cases}.$$ In addition, if the type of boundary condition in the problem is changed, the boundary conditions that Green's function needed to satisfy must also be changed to the corresponding zero boundary condition.
\end{definition}

By the definition of fundamental solution and Green's function, the differences between them is clear and easy to understand. Fundamental solution is only depend on the operator while Green's function is depend both on the operator and boundary.

\section{Operator Learning for PDE Systems} \label{Operator Learning for PDE Systems}
In this section, we consider the PDE system
\begin{equation} \label{eq0310_1}
    \begin{cases} u_1 - \lambda \Delta u_1 = f_1, \\ 
    \lambda \Delta u_2 + u_2 = f_2, \end{cases}
\end{equation}
where $u_1, u_2$ are unknown functions, $f_1, f_2$ are known source terms and $\lambda > 0$ is a parameter. In matrix-vector form, the system can be written as 
\begin{equation} \label{eq0310_2}
\mathbf{L}_\lambda \mathbf{u} = \mathbf{f},
\end{equation}
in which $\mathbf{u} = (u_1, u_2)^\top, \mathbf{f} = (f_1, f_2)^\top$ and 
\begin{equation} \label{eq0310_3}
\mathbf{L}_\lambda = \begin{pmatrix}
I & -\lambda \Delta \\
\lambda \Delta & I
\end{pmatrix}
\end{equation}
is a matrix of operators. Thus, a PDE system with Dirichlet boundary condition can be written as
\begin{equation} \label{eq0310_4}
    \begin{cases} \mathbf{L}_\lambda \mathbf{u} = \mathbf{f}, & \mathbf{x} \in \Omega \\ 
    \mathbf{u} = \mathbf{g}^D, & \mathbf{x} \in \partial \Omega, \end{cases}
\end{equation}
where $\mathbf{g}^D = (g_1, g_2)^\top$ is the vector of boundary conditions.

For system \eqref{eq0310_4}, we can also derive the corresponding boundary integral framework to solve it. The detailed derivation process is provided in the Appendix \ref{Boundary Integral Derivation for PDE System}. And we only present the results here. For operator $\mathbf{L}_\lambda$, the fundamental solution $\mathbf{G}_{\lambda,0}(\mathbf{x}, \mathbf{y}) = \{G^{i, j}_{\lambda,0}(\mathbf{x}, \mathbf{y})\}_{i, j = 1}^2$ is a $2 \times 2$ matrix such that 
\begin{equation} \label{eq0310_5}
{\mathbf{L}_\lambda}_x \mathbf{G}_{\lambda,0}(x, y) = \delta(x - y) \mathbf{Id}_2,
\end{equation}
where $\delta$ is the Dirac delta function and $\mathbf{Id}_2$ is the $2 \times 2$ identity matrix. And we can define the Green's function $\mathbf{G}_{\lambda}(\mathbf{x}, \mathbf{y})$  in a similar way as in Appendix \ref{appendix A}. A particular solution solution for $\mathbf{L}_\lambda\mathbf{u} = \mathbf{f}$ can be found as 
\begin{equation} \label{eq0310_6}
\mathcal{G}[\mathbf{f}](\mathbf{x}) := \int_\Omega \mathbf{G}_{\lambda}(\mathbf{x}, \mathbf{y}) \cdot \mathbf{f}(\mathbf{y}) d\mathbf{y},
\end{equation}
where `$\cdot$' represents the matrix-vector multiplication. In fact, $\mathcal{G}[\mathbf{f}]|_{\partial \Omega} \equiv \mathbf{0}$. The solution for \eqref{eq0310_4}
with $\mathbf{f} = \mathbf{0}$ can be given by
\begin{equation} \label{eq0310_7}
\mathcal{D}[\Phi] (\mathbf{x}) := - \int_{\partial \Omega} \frac{\partial \mathbf{K}_{\lambda}(\mathbf{x}, \mathbf{y})}{\partial n_{\mathbf{y}}} \cdot \Phi^*(\mathbf{y}) ds_{\mathbf{y}}, \quad \quad \mathbf{x} \in \Omega 
\end{equation}
where the vector $\Phi$ of density functions are the solution of this boundary integral equation
\begin{equation} \label{eq0310_8}
\frac{1}{2} \Phi(\mathbf{x}) - \int_{\partial \Omega} \frac{\partial \mathbf{K}_{\lambda}(\mathbf{x}, \mathbf{y})}{\partial n_{\mathbf{y}}} \cdot \Phi^*(\mathbf{y}) ds_{\mathbf{y}} = \mathbf{g}^D (\mathbf{y}), \quad \quad \mathbf{x} \in \partial \Omega.
\end{equation}
In the above statement, $\Phi = (\varphi_1, \varphi_2)^\top$ is a vector of density functions, $\Phi^* :=  (\varphi_2, \varphi_1)^\top$ and $\mathbf{K}_\lambda$ is defined by 
\begin{equation} \label{eq0310_9}
\mathbf{K}_\lambda = \begin{pmatrix}
\lambda G_{\lambda,0}^{1, 1} & -\lambda G_{\lambda,0}^{1, 2}\\ 
\lambda G_{\lambda,0}^{2, 1} & -\lambda G_{\lambda,0}^{2, 2} \end{pmatrix}.
\end{equation}
In summary, the solution to system \eqref{eq0310_4} can be expressed as the sum of $\mathcal{G}[\mathbf{f}]$ and $\mathcal{D}[\Phi]$, which correspond to the solutions of system \eqref{eq0310_4} under the constraints $\mathbf{g}^D = \mathbf{0}$ and $\mathbf{f} = \mathbf{0}$, respectively.

Building on the boundary integral theory discussed above, we can adopt an operator learning approach similar to that presented in Section \ref{Operator Learning}. Specifically, we construct two neural networks: $\mathcal{NN}_1$, which maps the source term $\mathbf{f}$ to the integral volume $\mathcal{G}[\mathbf{f}]$, and $\mathcal{NN}_2$, which maps the Dirichlet boundary condition $\mathbf{g}^D$ to the corresponding density function $\Phi$. As previously mentioned, the parameters of the PDE systems, specifically $\lambda$, can also be incorporated as inputs into these two networks.

\section{Derivation for Boundary Integral Theory to PDE Systems} \label{Boundary Integral Derivation for PDE System}
In this section, we present the detailed derivation of the boundary integral theory of operator $\mathbf{L}_\lambda$ of equation \eqref{eq0310_4} involved in Appendix \ref{Operator Learning for PDE Systems}.

\subsection{Derivation for The Fundamental Solution} \label{Derivation for The Fundamental Solution}
In the following derivation, we firstly fix $\lambda > 0$, and then briefly note that $G_{ij}(\mathbf{x}) = G^{i, j}_{\lambda, 0}(\mathbf{x}, \mathbf{0})$. From equation \eqref{eq0310_5}, we can obtain
\begin{equation} \label{eq0611_1}
\begin{cases}
G_{11}(\mathbf{x}) - \lambda \Delta G_{21}(\mathbf{x}) = \delta(\mathbf{x}) \\
G_{12}(\mathbf{x}) - \lambda \Delta G_{22}(\mathbf{x}) = 0 \\
\lambda \Delta G_{11}(\mathbf{x}) + G_{21}(\mathbf{x}) = 0 \\
\lambda \Delta G_{12}(\mathbf{x}) + G_{22}(\mathbf{x}) = \delta(\mathbf{x}).
\end{cases}
\end{equation}
Applying the Fourier transform gives:
\begin{equation} \label{eq0611_2}
\begin{cases}
\hat{G}_{11}(\xi) + \lambda |\xi|^2 \hat{G}_{21}(\xi) = 1, \\
\hat{G}_{12}(\xi) + \lambda |\xi|^2 \hat{G}_{22}(\xi) = 0, \\
- \lambda |\xi|^2 \hat{G}_{11}(\xi) + \hat{G}_{21}(\xi) = 0, \\
- \lambda |\xi|^2 \hat{G}_{12}(\xi) + \hat{G}_{22}(\xi) = 1.
\end{cases}
\end{equation}
Solving the system, we obtain:
\begin{equation} \label{eq0611_3}
\begin{aligned} 
\hat{G}_{11}(\xi) &= \frac{1}{1 + \lambda^2 |\xi|^4}, &
\hat{G}_{12}(\xi) &= -\frac{\lambda |\xi|^2}{1 + \lambda^2 |\xi|^4}, \\
\hat{G}_{21}(\xi) &= \frac{\lambda |\xi|^2}{1 + \lambda^2 |\xi|^4}, &
\hat{G}_{22}(\xi) &= \frac{1}{1 + \lambda^2 |\xi|^4}.
\end{aligned}
\end{equation}
Then the inverse Fourier transform of $\hat{G}_{11}(\xi)$ gives that
\begin{align}
G_{11}(\mathbf{x}) = \mathcal{F}^{-1} \left( \hat{G}_{11} \right)(\mathbf{x})
& = \frac{1}{(2\pi)^2} \int_{\mathbb{R}^2} e^{i \mathbf{x} \cdot \xi} \frac{1}{1 + \lambda^2 |\xi|^4} \, d\xi \label{eq:fourier1} \\
&= \frac{1}{(2\pi)^2} \int_0^{2\pi} \int_0^\infty e^{i r |\mathbf{x}| \cos \theta} \frac{r}{1 + \lambda^2 r^4} \, dr \, d\theta \label{eq:fourier2} \\
&= \frac{1}{2\pi} \int_0^\infty J_0(|\mathbf{x}| r) \frac{r}{1 + \lambda^2 r^4} \, dr \label{eq:fourier3} \\
&= -\frac{1}{2\pi \lambda} \textbf{kei}_0 \left( \frac{|\mathbf{x}|}{\sqrt{\lambda}} \right), \label{eq:fourier4}
\end{align}
where $\textbf{kei}_0 (\textbf{x})$ is the Kelvin function defined by $ e^{-i \gamma \frac{\pi}{2}} K_\gamma \left( z e^{i \frac{\pi}{4}} \right) = \textbf{ker}_\gamma(z) + i \textbf{kei}_\gamma(z)$, and here $K_\gamma$ is the Modified Bessel function of the second kind of order $\gamma$. Similarly, we can obtain that 
\begin{equation} \label{eq0611_4}
G_{21}(\mathbf{x}) = \frac{1}{2\pi \lambda} \textbf{ker}_0 (\frac{|\mathbf{x}|}{\sqrt{\lambda}} ) 
\end{equation}
To conclude, by considering the radial symmetry of the fundamental solution, we arrive at
\begin{equation} \label{eq0611_5}
\begin{aligned} 
G_{\lambda, 0}^{1, 1} (\mathbf{x}, \mathbf{y}) &= - \frac{1}{2\pi \lambda} \textbf{kei}_0 (\frac{|\mathbf{x} - \mathbf{y}|}{\sqrt{\lambda}}), &
G_{\lambda, 0}^{1, 2} (\mathbf{x}, \mathbf{y}) &= - \frac{1}{2\pi \lambda} \textbf{ker}_0(\frac{|\mathbf{x} - \mathbf{y}|}{\sqrt{\lambda}}), \\
G_{\lambda, 0}^{2, 1} (\mathbf{x}, \mathbf{y}) &= \frac{1}{2\pi \lambda} \textbf{ker}_0(\frac{|\mathbf{x} - \mathbf{y}|}{\sqrt{\lambda}}), &
G_{\lambda, 0}^{2, 2} (\mathbf{x}, \mathbf{y}) &= - \frac{1}{2\pi \lambda} \textbf{kei}_0(\frac{|\mathbf{x }- \mathbf{y}|}{\sqrt{\lambda}}).
\end{aligned}
\end{equation}

\subsection{Derivation for The Boundary Integral Equation} \label{Derivation for The Boundary Integral Equation}
In this section, we want to obtain the boundary integral equation of PDE system \eqref{eq0310_4} with $\mathbf{f} = \mathbf{0}$. For $u, v, p, q \in C^2(\bar{\Omega}; \mathbb{R})$:
\begin{equation} \label{eq0611_6}
\int_\Omega (p, q) \cdot \mathcal{L} \begin{pmatrix} u \\ v \end{pmatrix} = \int_\Omega \left( up + vq \right) + \lambda \int_\Omega \left( q \Delta u - p \Delta v \right), \\
\end{equation}
\begin{equation} \label{eq0611_7}
\int_\Omega (u, -v) \cdot \mathcal{L} \begin{pmatrix} p \\ -q \end{pmatrix} = \int_\Omega (up + vq)  + \lambda \int_\Omega ( - v \Delta p + u \Delta q ), \\
\end{equation}
By equation \eqref{eq0611_6} - \eqref{eq0611_7} and Green’s second identity, we have
\begin{equation} \label{eq0611_8}
\int_\Omega (p, q) \cdot \mathcal{L} \begin{pmatrix} u \\ v \end{pmatrix} - (u, -v) \cdot \mathcal{L} \begin{pmatrix} p \\ -q \end{pmatrix} d\mathbf{x} = \lambda \int_{\partial \Omega} (\frac{\partial u}{\partial n}q - \frac{\partial q}{\partial n}u) ds + \lambda \int_{\partial \Omega} (\frac{\partial p}{\partial n}v - \frac{\partial v}{\partial n}p) ds.
\end{equation}
Now, let $p(\mathbf{x}) = G^{1,1}_{\lambda, 0}(\mathbf{x}, \mathbf{y}), q(\mathbf{x}) = G^{1, 2}_{\lambda, 0}(\mathbf{x}, \mathbf{y})$ in \eqref{eq0611_8}, consider the property of $\delta$-function, we have
\begin{align} \label{eq0611_9}
&\quad \int_\Omega  
\begin{pmatrix}
G^{1,1}_{\lambda, 0}(\mathbf{x}, \mathbf{y}) \\
G^{1,2}_{\lambda, 0}(\mathbf{x}, \mathbf{y})
\end{pmatrix}^\top \cdot
\mathcal{L} 
\begin{pmatrix} 
u \\ v 
\end{pmatrix} 
(\mathbf{x}) 
\, d\mathbf{x} \notag \\
&= \lambda \int_{\partial \Omega} 
\left[
\frac{\partial u (\mathbf{x})}{\partial n_{\mathbf{x}}} G^{1,2}_{\lambda, 0} (\mathbf{x}, \mathbf{y})
- \frac{\partial G^{1,2}_{\lambda, 0} (\mathbf{x}, \mathbf{y})}{\partial n_\mathbf{x}} u(\mathbf{x})
+ \frac{\partial G^{1,1}_{\lambda, 0} (\mathbf{x}, \mathbf{y})}{\partial n_\mathbf{x}} v(\mathbf{x})
- \frac{\partial v(\mathbf{x})}{\partial n_{\mathbf{x}}} G^{1,1}_{\lambda, 0}(\mathbf{x}, \mathbf{y})
\right] ds \notag \\
&+ 
\begin{cases} 
u(\mathbf{y}), & \mathbf{y} \in \Omega, \\
\frac{1}{2} u(\mathbf{y}), & \mathbf{y} \in \partial \Omega, \\
0, & \mathbf{y} \in \bar{\Omega}^c.
\end{cases}
\end{align}
Similarly, let $p(\mathbf{x}) = G^{1,2}_{\lambda, 0}(\mathbf{x}, \mathbf{y}), q(\mathbf{x}) = -G^{2, 2}_{\lambda, 0}(\mathbf{x}, \mathbf{y})$ in \eqref{eq0611_8}, and consider $G^{1,2}_{\lambda, 0} = - G^{2, 1}_{\lambda, 0}$, we have
\begin{align} \label{eq0611_10}
&\quad\int_\Omega  
\begin{pmatrix}
G^{2,1}_{\lambda, 0}(\mathbf{x}, \mathbf{y}) \\
G^{2,2}_{\lambda, 0}(\mathbf{x}, \mathbf{y})
\end{pmatrix}^\top \cdot
\mathcal{L} 
\begin{pmatrix} 
u \\ v 
\end{pmatrix} 
(\mathbf{x}) 
\, d\mathbf{x} \notag \\
&= \lambda \int_{\partial \Omega} 
\left[
\frac{\partial u (\mathbf{x})}{\partial n_\mathbf{x}} G^{2,2}_{\lambda, 0} (\mathbf{x}, \mathbf{y})
- \frac{\partial G^{2,2}_{\lambda, 0} (\mathbf{x}, \mathbf{y})}{\partial n_\mathbf{x}} u(\mathbf{x})
+ \frac{\partial G^{2,1}_{\lambda, 0} (\mathbf{x}, \mathbf{y})}{\partial n_\mathbf{x}} v(\mathbf{x})
- \frac{\partial v(\mathbf{x})}{\partial n_\mathbf{x}} G^{2,1}_{\lambda, 0}(\mathbf{x}, \mathbf{y})
\right] ds \notag \\
&+ 
\begin{cases} 
v(\mathbf{y}), & \mathbf{y} \in \Omega, \\
\frac{1}{2} v(\mathbf{y}), & \mathbf{y} \in \partial \Omega, \\
0, & \mathbf{y} \in \bar{\Omega}^c.
\end{cases}
\end{align}
Combining \eqref{eq0611_9} and \eqref{eq0611_10}, we have
\begin{align} \label{eq0611_11}
&\quad \int_\Omega  
\mathbf{G}_{\lambda, 0}(\mathbf{x}, \mathbf{y}) \cdot
\mathcal{L} 
\begin{pmatrix} 
u \\ v 
\end{pmatrix} 
(\mathbf{x}) 
\, d\mathbf{x} \notag \\
& = \lambda \int_{\partial \Omega}   
\begin{pmatrix}
\frac{\partial u (\mathbf{x})}{\partial n_\mathbf{x}} G^{1,2}_{\lambda, 0} (\mathbf{x}, \mathbf{y})
- \frac{\partial G^{1,2}_{\lambda, 0} (\mathbf{x}, \mathbf{y})}{\partial n_\mathbf{x}} u(\mathbf{x})
+ \frac{\partial G^{1,1}_{\lambda, 0} (\mathbf{x}, \mathbf{y})}{\partial n_\mathbf{x}} v(\mathbf{x})
- \frac{\partial v(\mathbf{x})}{\partial n_\mathbf{x}} G^{1,1}_{\lambda, 0}(\mathbf{x}, \mathbf{y})
\\
\frac{\partial u (\mathbf{x})}{\partial n_\mathbf{x}} G^{2,2}_{\lambda, 0} (\mathbf{x}, \mathbf{y})
- \frac{\partial G^{2,2}_{\lambda, 0} (\mathbf{x}, \mathbf{y})}{\partial n_\mathbf{x}} u(\mathbf{x})
+ \frac{\partial G^{2,1}_{\lambda, 0} (\mathbf{x}, \mathbf{y})}{\partial n_\mathbf{x}} v(\mathbf{x})
- \frac{\partial v(\mathbf{x})}{\partial n_\mathbf{x}} G^{2,1}_{\lambda, 0}(\mathbf{x}, \mathbf{y})
\end{pmatrix}
ds \\
& + \begin{cases} 
\begin{pmatrix} u(\mathbf{y}) \\ v(\mathbf{y}) \end{pmatrix} , & \mathbf{y} \in \Omega, \\
\frac{1}{2} \begin{pmatrix} u(\mathbf{y}) \\ v(\mathbf{y}) \end{pmatrix}, & \mathbf{y} \in \partial \Omega, \\
\begin{pmatrix} 0 \\ 0 \end{pmatrix}, & \mathbf{y} \in \bar{\Omega}^c.
\end{cases}
\end{align}
Now, if $\mathbf{L}_\lambda (u, v)^\top = (0, 0)^\top$ in $\Omega$, then 
\begin{equation} \label{eq0611_13}
- \lambda \int_{\partial \Omega}   
\frac{\partial}{\partial n_\mathbf{x}} \widetilde{\mathbf{G}} (\mathbf{x}, \mathbf{y}) \cdot \begin{pmatrix} v \\ u \end{pmatrix} (\mathbf{x}) 
- \widetilde{\mathbf{G}}(\mathbf{x}, \mathbf{y}) \frac{\partial}{\partial n_\mathbf{x}} \begin{pmatrix} v \\ u \end{pmatrix} (\mathbf{x}) ds 
= 
\begin{cases} 
\begin{pmatrix} u(\mathbf{y}) \\ v(\mathbf{y}) \end{pmatrix} , & \mathbf{y} \in \Omega, \\
\frac{1}{2} \begin{pmatrix} u(\mathbf{y}) \\ v(\mathbf{y}) \end{pmatrix}, & \mathbf{y} \in \partial \Omega, \\
\begin{pmatrix} 0 \\ 0 \end{pmatrix}, & \mathbf{y} \in \bar{\Omega}^c,
\end{cases}
\end{equation}
where
\begin{equation} \label{eq0611_14}
\widetilde{\mathbf{G}} \triangleq \begin{pmatrix}
G^{1,1}_{\lambda, 0} \quad -G^{1,2}_{\lambda, 0}  \\
G^{2,1}_{\lambda, 0} \quad -G^{2,2}_{\lambda, 0}
\end{pmatrix}.
\end{equation}
Similarly, if $\mathbf{L}_\lambda (p, q)^\top = (0, 0)^\top$ in $\bar{\Omega}^C$ and $p, q$ vanish at infinity,
\begin{equation} \label{eq0611_15}
- \lambda \int_{\partial \Omega}   
-\frac{\partial}{\partial n_\mathbf{x}} \widetilde{\mathbf{G}} (\mathbf{x}, \mathbf{y}) \cdot \begin{pmatrix} q \\ p \end{pmatrix} (\mathbf{x}) 
+ \widetilde{\mathbf{G}}(\mathbf{x}, \mathbf{y}) \frac{\partial}{\partial n_\mathbf{x}} \begin{pmatrix} q \\ p \end{pmatrix} (\mathbf{x}) ds 
= 
\begin{cases} 
\begin{pmatrix} 0 \\ 0 \end{pmatrix} , & \mathbf{y} \in \Omega, \\
\frac{1}{2} \begin{pmatrix} q(\mathbf{y}) \\ p(\mathbf{y}) \end{pmatrix}, & \mathbf{y} \in \partial \Omega, \\
\begin{pmatrix} q(\mathbf{y}) \\ p(\mathbf{y}) \end{pmatrix}, & \mathbf{y} \in \bar{\Omega}^c,
\end{cases}
\end{equation}
Denote $\begin{pmatrix} \varphi_1 \\ \varphi_2 \end{pmatrix} = \begin{pmatrix} u \\ v \end{pmatrix} - \begin{pmatrix} p \\ q \end{pmatrix}$ and $\begin{pmatrix} \psi_1 \\ \psi_2 \end{pmatrix} = \frac{\partial}{\partial n} \begin{pmatrix} u \\ v \end{pmatrix} - \frac{\partial}{\partial n} \begin{pmatrix} p \\ q \end{pmatrix}$, then \eqref{eq0611_14} and \eqref{eq0611_15} lead to
\begin{equation} \label{eq0611_16}
- \lambda \int_{\partial \Omega}   
\frac{\partial}{\partial n} \widetilde{\mathbf{G}} (\mathbf{x}, \mathbf{y}) \cdot \begin{pmatrix} \varphi_2 \\ \varphi_1 \end{pmatrix} (\mathbf{x}) 
- \widetilde{\mathbf{G}}(\mathbf{x}, \mathbf{y}) \frac{\partial}{\partial n} \begin{pmatrix} \psi_2 \\ \psi_1 \end{pmatrix} (\mathbf{x}) ds 
= 
\begin{cases} 
\begin{pmatrix} u(\mathbf{y}) \\ v(\mathbf{y}) \end{pmatrix} , & \mathbf{y} \in \Omega, \\
\frac{1}{2} (\begin{pmatrix} u(\mathbf{y}) \\ v(\mathbf{y}) \end{pmatrix} + \begin{pmatrix} q(\mathbf{y}) \\ p(\mathbf{y}) \end{pmatrix}), & \mathbf{y} \in \partial \Omega, \\
\begin{pmatrix} q(\mathbf{y}) \\ p(\mathbf{y}) \end{pmatrix}, & \mathbf{y} \in \bar{\Omega}^c,
\end{cases}
\end{equation}
As a matter of fact, supposing that $u, v$ are given, we can choose $p, q$ to match $u$ such that $\psi_1 = \psi_2 \equiv 0$, and then 
\begin{equation} \label{eq0611_17}
- \lambda \int_{\partial \Omega}   
\frac{\partial}{\partial n} \widetilde{\mathbf{G}} (\mathbf{x}, \mathbf{y}) \cdot \begin{pmatrix} \varphi_2 \\ \varphi_1 \end{pmatrix} (\mathbf{x}) ds 
= 
\begin{cases} 
\begin{pmatrix} u(\mathbf{y}) \\ v(\mathbf{y}) \end{pmatrix} , & \mathbf{y} \in \Omega, \\
\begin{pmatrix} u(\mathbf{y}) \\ v(\mathbf{y}) \end{pmatrix} - \frac{1}{2} \begin{pmatrix} \varphi_1(\mathbf{y}) \\ \varphi_2(\mathbf{y}) \end{pmatrix}, & \mathbf{y} \in \partial \Omega. \\
\end{cases}
\end{equation}
From this, we can get the boundary integral equation \eqref{eq0310_8}.

\end{document}